\documentclass[11pt]{article}

\newtheorem{lem}{Lemma}
\newtheorem{theo}{Theorem}
\newtheorem{rem}{Remark}
\newtheorem{prop}{Proposition}

\newenvironment{pf}{\ \\ {\bf Proof: }}{\hfill\mbox{$\diamond$}\medskip}
\usepackage{amsmath}
\usepackage{amssymb}

\begin{document}
\author{Vlad Bally$^*$
\and  Emmanuelle Cl\'ement\thanks{Laboratoire d'Analyse et de
Math\'ematiques Appliqu\'ees, UMR $8050$,
Universit\'e Paris-Est Marne-la-Vall\'ee, 5  Bld Descartes,
Champs-sur-marne, 77454 Marne-la-Vall\'ee Cedex 2, France. }}
\title{Integration by parts formula and applications to equations with jumps}
\date{Preliminary version, November 2009}
\maketitle

\begin{abstract}
We establish an integration by parts formula in an abstract framework in order to study the regularity 
of the law for processes solution of stochastic differential equations with jumps, including equations with discontinuous coefficients for which the Malliavin calculus
developed by                     
 Bismut and  Bichteler, Gravereaux and Jacod fails.

\textbf{2000 MSC}. Primary: 60H07, Secondary 60G51

\textbf{Key words}: Integration by parts formula, Malliavin Calculus,
Stochastic Equations, Poisson Point Measures.
\end{abstract}

\section{Introduction}
This paper is made up of two parts. In a first part we give an abstract, finite dimensional
version of Malliavin calculus. Of course Malliavin calculus is known to be
an infinite dimensional differential calculus and so  a finite dimensional
version seems to be of a limited interest. We discuss later on the relation
between the finite dimensional and the infinite dimensional frameworks and we highlight the
interest of the finite dimensional approach.

In the second part of the paper we use the results from the first section in
order to give sufficient conditions for the regularity of the law of $X_{t},$
where $X$ is the Markov process with infinitesimal operator
\begin{equation}
Lf(x)=\left\langle \nabla f(x),g(x)\right\rangle
+\int_{R^{d}}(f(x+c(z,x))-f(x))\gamma (z,x)\mu (dz).  \label{i0}
\end{equation}
Suppose for the moment that $\gamma $ does not depend on $x.$ Then it is well
known that the process $X$ may be represented as the solution of a
stochastic equation driven by a Poisson point measure with intensity measure 
$\gamma (z)\mu (dz).$ Sufficient conditions for the regularity of the law of 
$X_{t}$ using a Malliavin calculus for Poisson point measures are given in $
[B.G.J].$ But in our framework $\gamma $ depends on $x$ which roughly speaking
means that the law of the jumps depends on the position of the particle
 when the jump occurs. Such processes are of interest in a lot of
applications and unfortunately the standard Malliavin calculus developed in $[B.G.J]$ does
not apply in this framework. After the classical paper of Bichteler Gravereaux
and Jacod a huge work concerning the Malliavin calculus for Poisson point
measures has been done and many different approaches have been developed.
But as long as we know they do not lead to a solution for our problem. If $X$
is an one dimensional process an analytical argument permits to solve the
above problem , this is done in $[F.1],[F.2]$ and $[F.G]$ but the argument
there seems difficult to extend in the multi-dimensional case.

We come now back to the relation between the finite dimensional and the
infinite dimensional framework. This seems to be the more interesting point
in our approach so we try to explain the main idea. In order to prove
Malliavin's regularity criterion for the law of a functional $F$ on the
Wiener space the main tool is the integration by parts formula 
\begin{equation}
E(\partial ^{\beta }f(F))=E(f(F)H_{\beta })  \label{i1}
\end{equation}
where $\partial ^{\beta }$ denotes the derivative corresponding to a multi-index $\beta $ and $H_{\beta }$ is a random variable built using the
Malliavin derivatives of $F.$ Once such a formula is proved one may estimate
the Fourier transform $\widehat{p}_{F}(\xi )=E(\exp (i\xi F))$ in the
following way. First we remark that $\partial _{x}^{\beta }\exp (i\xi
x)=(i\xi )^{\beta }\exp (i\xi x)$ (with an obvious abuse of notation) and
then, using the integration by parts formula 
\begin{eqnarray*}
\left\vert \widehat{p}_{F}(\xi )\right\vert &=&\frac{1}{\vert \xi \vert^{\vert \beta \vert}}
\left\vert E(\partial _{x}^{\beta }\exp (i\xi F))\right\vert \\
&=&\frac{1}{\vert \xi \vert^{\vert \beta \vert}}\left\vert E(\exp (i\xi F)H_{\beta
})\right\vert \leq \frac{1}{\vert \xi \vert^{\vert \beta \vert}}E\left\vert H_{\beta
}\right\vert .
\end{eqnarray*}
If we know that $E\left\vert H_{\beta }\right\vert <\infty $ for every multi-index $\beta $ then we have proved that $\left\vert \xi \right\vert
^{p}\left\vert \widehat{p}_{F}(\xi )\right\vert $ is integrable for every $
p\in \mathbb{N}$ and consequently the law of $F$ is absolutely continuous with
respect to the Lebesgue measure and has an infinitely differentiable density.

Let us come back to the infinite dimensional differential calculus which
permits to built $H_{\beta }.$ In order to define the Malliavin derivative
of $F$ one considers a sequence of simple functionals $F_{n}\rightarrow F$
in $L^{2}$ and, if $DF_{n}\rightarrow G$ in $L^{2},$ then one defines $DF=G.$
The simple functionals $F_{n}$ are functions of a finite number of random
variables (increments of the Brownian motion) and the derivative $DF_{n}$ is
a gradient type operator defined in an elementary way. Then one may take the
following alternative way in order to prove the regularity of the law of $F.$
For each fixed $n$ one proves the analogues of the integration by parts
formula (\ref{i1}): $E(\partial ^{\beta }f(F_{n}))=E(f(F_{n})H_{\beta
}^{n}). $ As $F_{n}$ is a function which depends on a finite number $m$ of
random variables, such a formula is obtained using standard integration by
parts on $\mathbb{R}^{m}$ (this is done in the first section of this paper$).$ Then
the same calculus as above gives $\left\vert \widehat{p}_{F_{n}}(\xi
)\right\vert \leq \vert \xi \vert^{-\vert \beta \vert}E\left\vert H_{\beta }^{n}\right\vert .$
Passing to the limit one obtains 
\begin{equation*}
\left\vert \widehat{p}_{F}(\xi )\right\vert =\lim_{n}\left\vert \widehat{p}
_{F_{n}}(\xi )\right\vert \leq \vert \xi \vert^{-\vert \beta \vert} \sup_{n}E\left\vert H_{\beta
}^{n}\right\vert
\end{equation*}
and, if we can prove that $\sup_{n}E\left\vert H_{\beta }^{n}\right\vert
<\infty ,$ we are done. Notice that here we do not need that $
F_{n}\rightarrow F$ in $L^{2}$ but only in law. And also, we do not need to
built $H_{\beta }$ but only to prove that $\sup_{n}E\left\vert H_{\beta
}^{n}\right\vert <\infty .$ Anyway we are not very far from the standard
Malliavin calculus. Things become different if $\sup_{n}E\left\vert H_{\beta
}^{n}\right\vert =\infty $ and this is the case in our examples (because the
Ornstein Uhlenbeck operators $LF_{n}$ blow up as $n\rightarrow \infty ).$
But even in this case one may obtain estimates of the Fourier transform of $
F $ in the following way. One writes 
\begin{equation*}
\left\vert \widehat{p}_{F}(\xi )\right\vert \leq \left\vert \widehat{p}
_{F}(\xi )-\widehat{p}_{F_{n}}(\xi )\right\vert +\left\vert \widehat{p}
_{F_{n}}(\xi )\right\vert \leq \left\vert \xi \right\vert \times E\left\vert
F-F_{n}\right\vert +\vert \xi \vert^{-\vert \beta \vert}E\left\vert H_{\beta }^{n}\right\vert .
\end{equation*}
And if one may obtain a good balance between the convergence to zero of  the error $E\left\vert
F-F_{n}\right\vert $ and the blow up to infinity of $E\left\vert H_{\beta
}^{n}\right\vert  $ then one obtains $\left\vert \widehat{p
}_{F}(\xi )\right\vert \leq \left\vert \xi \right\vert ^{-p}$ for some $p$.
Examples in which such a balance works are given in Section 3.
An other application of this methodology is given in [B.F] for the Boltzmann equation. In this case some specific
and nontrivial difficulties appear due to the singularity and unboundedness of the coefficients of the equation.

The paper is organized as follows. In Section 2 we establish the abstract
Malliavin calculus associated to a finite dimensional random variable and we
obtain estimates of the weight $H_{\beta }$ which appear in the integration
by parts formula (we follow here some ideas which already appear in $
[B],[B.B.M]$ and $[Ba.M])$. Section 3 is devoted to the study of the regularity of the law of the Markov process $X$ of
infinitesimal operator (\ref{i0}) and it contains our main results : Proposition \ref{fourierX} and Theorem \ref{density}.
At last we provide in Section 4 the technical
estimates which are needed to prove the results of section 3.

\section{Integration by parts formula}
\subsection{Notations-derivative operators}

Throughout this paper, we consider a sequence of random variables $(V_{i})_{i\in \mathbb{N}^*}$ on a probability space $(\Omega ,\mathcal{F},P)$,
a sub $\sigma$-algebra $\mathcal{G}\subseteq \mathcal{F}$ and a random variable $J$, $\mathcal{G}$ measurable, with values in $\mathbb{N}$. We assume that
the variables $(V_i)$ and $J$ satisfy the following integrability conditions :
for all $p \geq 1$, $E(J^{p})+E((\sum_{i=1}^{J}V_{i}^{2})^{p})<\infty .$ 
Our aim is to establish a differential calculus based on the variables $(V_i)$, conditionally on $\mathcal{G}$, and we first define the class of functions on which this differential calculus
will apply. More precisely, we consider in this paper functions $f :\Omega \times \mathbb{R}^{\mathbb{N}^*}\rightarrow \mathbb{R}$  which can be written as 
\begin{equation}
f(\omega ,v)=\sum_{j=1}^{\infty }f^{j}(\omega
,v_{1},...,v_{j})1_{\{J(\omega )=j\}} \label{defS}
\end{equation}
where $f^{j}:\Omega \times \mathbb{R}^{j}\rightarrow \mathbb{R}$ are $\mathcal{G\times B}
(\mathbb{R}^{j})\mathcal{-}$measurable functions. We denote by $\mathcal{M}$ the class of  functions $
f$ given by (\ref{defS}) such that there exists a random variable $C\in \cap _{q\geq1}L^{q}(\Omega ,
\mathcal{G},P)$ and a real number $p\geq1$ satisfying $\left\vert f(\omega ,v)\right\vert \leq C(\omega
)(1+(\sum_{i=1}^{J(\omega )}v_{i}^{2})^{p})$. So conditionally on $\mathcal{G}$, the functions of $\mathcal{M}$
have polynomial growth with  respect to the variables $(V_i)$. We need some more
notations.  Let $\mathcal{G}_i$ be the $\sigma-$algebra generated by $\mathcal{G} \cup \sigma(V_j, 1 \leq j \leq J, j\neq i)$ and let $(a_i(\omega))$ and 
$(b_i(\omega))$ be
sequences of  $\mathcal{G}_i$ measurable random variables satisfying $-\infty \leq a_i(\omega)<b_i(\omega) \leq + \infty$, for all $i \in \mathbb{N}^*$. Now let $O_i$ be the open set of $\mathbb{R}^{\mathbb{N}^*}$ defined by $O_i=P_i^{-1}(]a_i,b_i[)$, where $P_i$ is the coordinate  map $P_i(v)=v_i$. We localize the differential calculus on the sets $(O_i)$ by introducing some weights $(\pi_i)$, satisfying the following hypothesis.

\textbf{H0.} For all $i \in \mathbb{N^*}$, $\pi_i \in \mathcal{M}$,  $0\leq \pi _{i}\leq 1$ and $\{\pi_i>0 \} \subset O_i$. Moreover for all $j \geq 1$, $\pi_i^j$ is infinitely differentiable with bounded
derivatives with respect to the variables $(v_1, \ldots, v_j)$. 

We associate to these weights $(\pi_i)$,  the spaces 
 $C_{\pi }^{k}\subset \mathcal{M},k\in \mathbb{N}^*$
defined recursively as follows. For $k=1$, $C_{\pi }^{1}$ denotes the space of  functions $f\in 
\mathcal{M}$ such that for each $i\in \mathbb{N}^*,$ $f$ admits a partial derivative with respect to the variable
$v_i$ on the open set $O_i$. We then define
\begin{equation*}
\partial _{i}^{\pi }f(\omega ,v):=\pi _{i}(\omega,v)
\frac{\partial 
}{\partial v_{i}}f(\omega ,v)
\end{equation*}
and we assume that $\partial _{i}^{\pi }f \in \mathcal{M}.$

Note that the chain rule is verified : for each $\phi \in C^{1}(\mathbb{R}^{d},\mathbb{R})$ and $
f=(f^{1},...,f^{d})\in (C_{\pi }^{1})^{d}$ we have 
\begin{equation*}
\partial _{i}^{\pi }\phi (f)=\sum_{r=1}^{d}\partial _{r}\phi
(f)\partial _{i}^{\pi }f^{r}.
\end{equation*}
Suppose now that $C_{\pi }^{k}$ is already defined. For a multi-index $
\alpha =(\alpha _{1},...,\alpha _{k})\in \mathbb{N}^{*k}$ 
we define recursively $\partial _{\alpha }^{\pi
}=\partial _{\alpha _{k}}^{\pi }...\partial _{\alpha _{1}}^{\pi }$ and
$C_{\pi }^{k+1}$ is the space of functions $f\in C_{\pi
}^{k}$ such that for every multi-index $\alpha =(\alpha _{1},...,\alpha
_{k})\in \mathbb{N}^{*k}$ we have $\partial _{\alpha }^{\pi }f\in C_{\pi }^{1}.$
Note that if $f\in C_{\pi
}^{k}$, $\partial _{\alpha }^{\pi }f \in \mathcal{M}$ for each $\alpha$ with $\left\vert \alpha
\right\vert \leq k$.

 Finally we define $C_{\pi }^{\infty
}=\cap _{k\in \mathbb{N}^*}C_{\pi }^{k}.$ Roughly speaking the space $C_{\pi }^{\infty
} $ is the analogue of $C^{\infty }$ with partial derivatives $\partial
_{i}$  replaced by  localized derivatives $\partial _{i}^{\pi }.$

\textbf{Simple functionals}. A random variable $F$ is called a simple
functional if there exists $f \in C_{\pi }^{\infty }$ such that $
F=f(\omega ,V)$, where $V=(V_i)$. We denote by $\mathcal{S}$ the space of the simple
functionals. Notice that $\mathcal{S}$ is an algebra. It is worth to remark that conditionally on $\mathcal{G}$, $F=f^J(V_1, \ldots, V_J)$.

\textbf{Simple processes}. A simple process is a sequence of random variables $
U=(U_{i})_{i\in \mathbb{N}^*}$ such that for each $i\in \mathbb{N}^*,$ $U_{i} \in \mathcal{S}$. 
Consequently, conditionally on $\mathcal{G}$, we have $U_i=u_i^J(V_1, \ldots, V_J)$.
We denote by $\mathcal{P}$\ the space of the
simple processes and we define the scalar product

\begin{equation*}
\left\langle U,V\right\rangle_J =\sum_{i=1}^{J}U_{i}V_{i}.
\end{equation*}
Note that $\left\langle U,V\right\rangle_J  \in \mathcal{S}$.

We can now  define the derivative operator and state the integration by parts formula.

$\square $ \textbf{The derivative operator.} We define  
$D:\mathcal{S}\rightarrow 
\mathcal{P}:$ by
\begin{equation*}
DF:=(D_{i}F)\in 
\mathcal{P} \quad  \mbox{where } D_{i}F:=\partial _{i}^{\pi }f(\omega ,V).       
\end{equation*}
Note that $D_iF=0$ for $i>J$.
For $F=(F^{1},...,F^{d})\in \mathcal{S}^{d}$ the Malliavin covariance matrix
is defined by
\begin{equation*}
\sigma ^{k,k'}(F)=\left\langle DF^{k},DF^{k'}\right\rangle_J =\sum_{j=1}^{J} D_{j}F^{k}D_{j}F^{k'}.
\end{equation*}
 We denote
\begin{equation*}
\Lambda (F)=\{\det \sigma (F)\neq 0\}\quad and\quad \gamma
(F)(\omega )=\sigma ^{-1}(F)(\omega ),\omega \in \Lambda
(F).
\end{equation*}

In order to derive an integration by parts formula, we need some additional assumptions on the random variables $(V_i)$.
The main hypothesis  is that conditionally on $\mathcal{G},$  the law of the vector$
(V_{1},...,V_{J})$ admits a locally smooth density with respect to the Lebesgue measure on $\mathbb{R}^J$.

\textbf{H1.} i) Conditionally on $\mathcal{G}$,  the vector $
(V_{1},...,V_{J})$ is absolutely continuous with respect to the Lebesgue measure on $\mathbb{R}^J$ and we note $p_J$ the conditional density.

\hspace{1cm} ii)  The set $\{p_J>0 \}$ is  open in $\mathbb{R}^J$ and on  $\{p_J>0 \}$ $\ln p_{J}\in C_{\pi }^{\infty }$.

\hspace{1cm} iii)  $\forall q \geq 1$, there exists a constant $C_q$ such that
$$
(1+ \vert v\vert^q)p_J \leq C_q
$$
where $\vert v\vert$ stands for the euclidian norm of the vector $(v_1, \ldots, v_J)$.

Assumption iii) implies in particular that
conditionally on $\mathcal{G},$ the functions of $\mathcal{M}$ are integrable with respect to $p_J$ and that for $f \in \mathcal{M}$ :
\begin{equation*}
E_{\mathcal{G}}(f(\omega ,V) )=\int_{\mathbb{R}^{J}}f^J \times
p_{J}(\omega ,v_{1},...,v_{J})dv_{1}...dv_{J}.
\end{equation*}

$\square $ \textbf{The divergence operator} \ Let $U=(U_{i})_{i\in \mathbb{N}^*}\in 
\mathcal{P}$ with $U_{i}\in \mathcal{S}.$ We define $
\delta :\mathcal{P}\rightarrow \mathcal{S}$ by
\begin{eqnarray*}
\delta _{i}(U) &:&=-(\partial _{v_i}(\pi _{i}U_{i})+U_{i} 1_{ \{p_J>0 \}} \partial_i^{\pi} \ln p_J), \\
\delta (U) &=&\sum_{i=1}^{J}\delta _{i}(U)
\end{eqnarray*}
For $F \in \mathcal{S}$, let $L(F)=\delta(DF)$. 
\subsection{Duality and integration by parts formulae}
In our framework the duality between $\delta $ and $D$ is given by the
following proposition.

\begin{prop}
Assume  \textbf{H0} and \textbf{H1}, then $\forall F\in 
\mathcal{S}$ and $\forall U\in \mathcal{P}$ we have 
\begin{equation}
E_{\mathcal{G}}(\left\langle DF,U\right\rangle_J )=E_{\mathcal{G}
}(F\delta (U)).  \label{1.2}
\end{equation}
\end{prop}

\begin{pf}
By definition, we have $E_{\mathcal{G}}(\left\langle DF,U\right\rangle_J )=\sum_{i=1}^{J} E_{\mathcal{G}}(D_i F \times U_i)$
and from \textbf{H1} 
\begin{equation*}
E_{\mathcal{G}}(D_{i}F\times U_{i})=\int_{\mathbb{R}^J} \partial
_{v_i}(f ^J) \pi_i \ u_{i}^J \ p_{J}(\omega
,v_{1},...,v_{J})dv_{1}...dv_{J}  
\end{equation*} 
 recalling that $\{\pi _{i}>0\}\subset O_i$, we obtain from Fubini's theorem
\begin{equation*}
 E_{\mathcal{G}}(D_{i}F\times U_{i})=\int_{\mathbb{R}^{J-1}} \left(\int_{a_i}^{b_i} \partial
_{v_i}(f ^J) \pi_i \ u_{i}^J \ p_{J}(\omega
,v_{1},...,v_{J})dv_i \right)dv_{1}..dv_{i-1}dv_{i+1}...dv_{J}.
\end{equation*}
By using the classical integration by parts formula, we have
$$
\int_{a_i}^{b_i} \partial
_{v_i}(f ^J) \pi_i \ u_{i} ^J\ p_{J}(\omega
,v_{1},...,v_{J})dv_i=[f^J \pi_i u_i^J p_J]_{a_i}^{b_i}-\int_{a_i}^{b_i}f^J \partial_{v_i}(u_i^J \pi_i p_J)dv_i.
$$
Now if $-\infty<a_i<b_i<+ \infty$, we have $\pi_i(a_i)=0=\pi_i(b_i)$ and $[f^J \pi_i u_i^J p_J]_{a_i}^{b_i}=0$. Moreover
since $f^J$, $u_i^J$ and $\pi_i$ belong to $\mathcal{M}$, we deduce from $H_1$ iii) that $\lim_{\vert v_i \vert \rightarrow +\infty} (f^J \pi_i u_i^J p_J)=0$ and we
obtain that for all $a_i$, $b_i$ such that $-\infty \leq a_i<b_i \leq+ \infty$ :
$$
\int_{a_i}^{b_i} \partial
_{v_i}(f ^J) \pi_i \ u_{i} ^J\ p_{J}(\omega
,v_{1},...,v_{J})dv_i=-\int_{a_i}^{b_i}f^J \partial_{v_i}(u_i^J \pi_i p_J)dv_i,
$$
Observing that $\partial_{v_i}(u_i^J \pi_i p_J)=(\partial_{v_i}(u_i^J \pi_i ) +u_i^J  1_{ \{p_J>0 \} }  \partial_i^{\pi}( \ln p_J ))p_J$,
 the proposition is proved. 
 
\end{pf}

We have the following straightforward computation rules.

\begin{lem}\label{CR}
Let $\phi :\mathbb{R}^{d}\rightarrow \mathbb{R}$ be a smooth function and $F=(F^{1},...,F^{d})
\in \mathcal{S}^{d}.$ Then $\phi (F)\in \mathcal{S}$ and 
\begin{equation}
D\phi (F)=\sum_{r=1}^{d}\partial _{r}\phi (F)DF^{r}.  \label{1.3}
\end{equation}
If $F\in \mathcal{S}$ and $U\in \mathcal{P}$ then
\begin{equation}
\delta (FU)=F\delta (U)-\left\langle DF,U\right\rangle_J .  \label{1.4}
\end{equation}
Moreover, for $F=(F^{1},...,F^{d})
\in \mathcal{S}^{d}$, we have
\begin{equation}
L \phi(F)=\sum_{r=1}^{d}\partial _{r}\phi (F)LF^{r} -\sum_{r,r'=1}^{d}\partial _{r,r'}\phi (F)\left\langle DF^{r}, DF^{r'}\right\rangle_J.  \label{CRL}
\end{equation}
\end{lem}

The first equality  is a consequence of the chain rule,  the 
second one follows from  the definition of the divergence operator $\delta$. Combining these equalities $(\ref{CRL})$ follows. 

We can now state the main results of this section.

\begin{theo} \label{IPP}
We assume  \textbf{H0} and \textbf{H1}. Let $F=(F^{1},...,F^{d})\in 
\mathcal{S}^{d}$, $G\in \mathcal{S}$ and $\phi :\mathbb{R}^{d}\rightarrow \mathbb{R}$ be a smooth  bounded function with
bounded derivatives. Let $
\Lambda \in \mathcal{G},\Lambda \subset \Lambda (F)$ such that 
\begin{equation}
E(\left\vert \det \gamma (F)\right\vert ^{p}1_{\Lambda })<\infty
\quad \forall p \geq 1.  \label{1.5}
\end{equation}
 Then, for every $r=1,...,d,$  
\begin{equation}
E_{\mathcal{G}}\left(\partial _{r}\phi (F)G\right)1_{\Lambda }=E_{\mathcal{G}}\left(\phi
(F)H_{r}(F,G)\right)1_{\Lambda }  \label{IPP1}
\end{equation}
with 
\begin{equation}
H_{r}(F,G)=\sum_{r'=1}^{d}\delta (G\gamma ^{r',r}(F)DF^{r'})
=\sum_{r'=1}^{d}\left(G \delta(\gamma ^{r',r}(F) DF^{r'})-\gamma ^{r',r} \left\langle
DF^{r'},D G\right\rangle_J  \right).  \label{weight1}
\end{equation}
\end{theo}

\begin{pf}Using the chain rule 
\begin{eqnarray*}
\left\langle D\phi (F), DF^{r'}\right\rangle_J  &=&\sum_{j=1}^{J} D_j \phi(F)
D_{j} F^{r'} \\
&=&\sum_{j=1}^{J}(\sum_{r=1}^{d}\partial _{r}\phi (F)D_j F^r)D_{j}F^{r'}
=\sum_{r=1}^{d}\partial _{r}\phi (F)\sigma ^{r,r'}(F)
\end{eqnarray*}
so that $\partial _{r}\phi (F)1_{\Lambda } = 1_{\Lambda } \sum_{r'=1}^{d}\left\langle D\phi (F),
DF^{r'}\right\rangle_J \gamma ^{r',r}(F).$ Since $F\in \mathcal{S}
^{d} $ it follows that $\phi (F)\in \mathcal{S}$ and $\ \sigma
^{r,r'}(F)\in \mathcal{S}.$ Moreover, since $\det \gamma (F)1_{\Lambda
}\in \cap _{p\geq 1}L^{p}$ it follows that $\gamma ^{r,r'}(F)1_{\Lambda
}\in \mathcal{S}.$ So $G\gamma ^{r',r}(F) DF^{r'}1_{\Lambda }\in 
\mathcal{P}$ and the duality formula gives: 
\begin{eqnarray*}
E_{\mathcal{G}}\left(\partial _{r}\phi (F)G\right)1_{\Lambda } &=&\sum_{r'=1}^{d}E_{
\mathcal{G}}\left(\left\langle D\phi (F),G\gamma ^{r',r}(F) DF^{r'}
\right\rangle_J \right)1_{\Lambda } \\
&=&\sum_{r'=1}^{d}E_{\mathcal{G}}\left(\phi (F)\delta (G\gamma ^{r',r}(F) DF^{r'}
)\right)1_{\Lambda }.
\end{eqnarray*}
\end{pf}

We can extend this integration by parts formula. 
\begin{theo} \label{IPPE}
Under the assumptions of Theorem \ref{IPP}, we have for every multi-index $\beta=(\beta_1, \ldots, \beta_q) \in \{1, \ldots, d\}^q$  
\begin{equation}
E_{\mathcal{G}}\left(\partial _{\beta}\phi (F)G\right)1_{\Lambda }=E_{\mathcal{G}}\left(\phi
(F)H_{\beta}^q(F,G)\right)1_{\Lambda }  \label{IPP2}
\end{equation}
where the weights $H^q$ are defined recursively by (\ref{weight1}) and
\begin{equation}
H_{\beta}^q(F,G)=H_{\beta_1}\left(F, H_{(\beta_2, \ldots,\beta_q)}^{q-1}(F,G) \right). \label{weight2}
\end{equation}
\end{theo}

\begin{pf}
The proof is straightforward by induction.
For $q=1$, this is just Theorem \ref{IPP}. Now assume that Theorem \ref{IPPE} is true for $q \geq 1$ and let us prove it for $q+1$. Let $\beta=(\beta_1, \ldots, \beta_{q+1}) \in \{1, \ldots, d\}^{q+1}$, we have
$$
E_{\mathcal{G}}\left(\partial _{\beta}\phi (F)G\right)1_{\Lambda }=
E_{\mathcal{G}}\left(\partial _{(\beta_2, \ldots, \beta_{q+1})}(\partial_{\beta_1}\phi (F))G\right)1_{\Lambda }
=E_{\mathcal{G}}\left(\partial _{\beta_1}\phi (F)H^q_{(\beta_2, \ldots, \beta_{q+1})}(F,G)\right)1_{\Lambda }
$$
and the result follows.
\end{pf}

\subsection{Estimations of $H^q$}
\subsubsection{Iterated derivative operators, Sobolev norms}
In order to estimate the weights $H^q$ appearing in the integration by parts formulae of the previous
section, we need first to define iterations of the derivative operator.
Let $\alpha=(\alpha_1, \ldots,\alpha_k)$ be a multi-index, with $\alpha_i \in \{1, \ldots,J\}$, for $i=1, \ldots,k$ and $\vert \alpha \vert=k$.
For $F \in \mathcal{S}$, we define recursively $D^k_{(\alpha_1, \ldots, \alpha_k)} F=D_{\alpha_k}
(D^{k-1}_{(\alpha_1, \ldots, \alpha_{k-1})} F)$ and 
$$D^k F=\left( D^k_{(\alpha_1, \ldots, \alpha_k)} F \right)_{\alpha_i \in \{1, \ldots,J\}}.$$
Remark that $D^k F \in \mathbb{R}^{J \otimes k}$ and consequently we define the norm of $D^k F$ as
$$
\vert D^k F\vert= \sqrt{\sum_{\alpha_1, \ldots, \alpha_k=1}^J \vert D^k_{(\alpha_1, \ldots, \alpha_k)}F \vert^2}.
$$
Moreover, we introduce the following norms, for $F \in \mathcal{S}$:
\begin{equation}
\vert F \vert_{1,l}=\sum_{k=1}^l \vert D^k F\vert, \hspace{0,5cm}
\vert F \vert_l= \vert F\vert + \vert F \vert_{1,l} = \sum_{k=0}^l \vert D^k F \vert .\label{norml}
\end{equation}
For $F=(F^1, \ldots , F^d) \in \mathcal{S}^d$:
$$
\vert F \vert_{1,l}=\sum_{r=1}^d \vert F^r \vert_{1,l}, \hspace{0,5cm} \vert F \vert_l=\sum_{r=1}^d \vert F^r \vert_l ,
$$
and similarly for $F=(F^{r,r'})_{r,r'=1, \ldots, d}$
$$
\vert F \vert_{1,l}=\sum_{r,r'=1}^d \vert F^{r,r'} \vert_{1,l} , \hspace{0,5cm} \vert F \vert_l=\sum_{r,r'=1}^d \vert F^{r,r'} \vert_l .
$$
Finally for $U=(U_i)_{i \leq J} \in \mathcal{P}$, we have $D^kU=(D^kU_i)_{i \leq J}$ and we define the norm of $D^kU$ as
$$
\vert D^kU \vert=\sqrt{\sum_{i=1}^J \vert D^kU_i \vert^2 }.
$$
We can remark that for $k=0$, this gives $\vert U \vert=\sqrt{\left\langle U,U \right\rangle_J }$. Similarly to (\ref{norml}), we set
\begin{equation*}
\vert U \vert_{1,l}=\sum_{k=1}^l \vert D^k U\vert, \hspace{0,5cm}
\vert U \vert_l= \vert U\vert + \vert U \vert_{1,l} = \sum_{k=0}^l \vert D^k U \vert .
\end{equation*}
Observe that for $F, G \in \mathcal{S}$, we have $D(F\times G)=DF \times G +F \times DG$. This leads to the following useful inequalities
\begin{lem}\label{ineq}
Let $F,G \in \mathcal{S}$ and $U,V \in \mathcal{P}$, we have
\begin{equation}
\vert F \times G \vert_l \leq 2^{l}\sum_{l_1+l_2\leq l} \vert F \vert_{l_1} \vert G \vert_{l_2}, \label{prod}
\end{equation}
\begin{equation}
\vert \left\langle U,V\right\rangle_J \vert_l \leq  2^{l} \sum_{l_1+l_2\leq l} \vert U \vert_{l_1} \vert V \vert_{l_2}. \label{scalP}
\end{equation}
\end{lem}
We can remark that the first inequality  is sharper than the following one $\vert F \times G \vert_l \leq C_l \vert F \vert_l \vert G \vert_l$.
Moreover from (\ref{scalP}) with $U=DF$ and $V=DG$ ( $F,G, \in \mathcal{S}$) we deduce
\begin{equation}
\vert \left\langle DF,DG \right\rangle_J \vert_l \leq  2^{l}\sum_{l_1+l_2\leq l} \vert F \vert_{1,l_1+1} \vert G \vert_{1,l_2+1} \label{scal1}
\end{equation}
and as an immediate consequence of (\ref{prod}) and (\ref{scal1}), we have for $F,G,H \in \mathcal{S}$:
\begin{equation}
\vert H \left\langle DF, DG \right\rangle_J \vert_l \leq 2^{2l}\sum_{l_1+l_2 +l_3\leq l} \vert F \vert_{1,l_1+1} \vert G \vert_{1,l_2+1} \vert H \vert_{l_3}. \label{scal3}
\end{equation}

\begin{pf}
We just prove (\ref{scalP}), since (\ref{prod}) can be proved on the same way.
We first give a bound for $D^k\left\langle U,V \right\rangle_J=(D^k_{\alpha}\left\langle U,V \right\rangle_J)_{ \alpha \in \{1,\ldots,J\}^k}$. For a multi-index $\alpha
=(\alpha _{1},...,\alpha _{k})$, with $\alpha_i \in \{1, \ldots, J\}$, we note $\alpha (\Gamma )=(\alpha
_{i})_{i\in \Gamma }$, where $\Gamma \subset \{1, \ldots, k \}$ and  $\alpha (\Gamma^c)=(\alpha
_{i})_{i\notin \Gamma }$. We have
\begin{equation*}
D_{\alpha }^{k}\left\langle U,V\right\rangle_J =\sum_{i=1}^{J}D_{\alpha
}^{k}(U_{i}V_{i})=\sum_{k'=0}^{k}\sum_{\left\vert \Gamma
\right\vert =k'} \sum_{i=1}^{J} D_{\alpha (\Gamma )}^{k'}U_i\times D_{\alpha (\Gamma^c)}^{k-k'}V_i.
\end{equation*}
Let $W^{i,\Gamma}=(W^{i ,\Gamma}_{\alpha })_{\alpha \in \{1,\ldots,J\}^k}=(D_{\alpha (\Gamma )}^{k'}U_i\times D_{\alpha (\Gamma^c)}^{k-k'}V_i)_{\alpha \in \{1,\ldots,J\}^k}$, we have the equality in $\mathbb{R}^{J \otimes k}$ :
$$
D^{k}\left\langle U,V\right\rangle_J =\sum_{k'=0}^{k}\sum_{\left\vert \Gamma
\right\vert =k'} \sum_{i=1}^{J} W^{i,\Gamma}.
$$
This gives
$$
\vert D^k \left\langle U,V\right\rangle_J \vert \leq \sum_{k'=0}^{k}\sum_{\left\vert \Gamma
\right\vert =k'} \vert \sum_{i=1}^{J} W^{i,\Gamma}\vert,
$$
where
$$
\vert \sum_{i=1}^{J} W^{i,\Gamma}\vert=\sqrt{\sum_{\alpha_1,\ldots,\alpha_k=1}^J \vert \sum_{i=1}^{J} W_{\alpha}^{i,\Gamma}\vert^2}.
$$
But from Cauchy Schwarz inequality, we have
$$
\vert \sum_{i=1}^{J} W_{\alpha}^{i,\Gamma}\vert^2=
\vert\sum_{i=1}^{J} D_{\alpha (\Gamma )}^{k'}U_i\times D_{\alpha (\Gamma^c)}^{k-k'}V_i\vert^2 \leq
\sum_{i=1}^{J} \vert D_{\alpha (\Gamma )}^{k'}U_i\vert^2 \times \sum_{i=1}^J \vert D_{\alpha (\Gamma^c)}^{k-k'}V_i \vert^2.
$$
Consequently we obtain
\begin{eqnarray*}
\vert \sum_{i=1}^{J} W^{i,\Gamma}\vert&\leq&
\sqrt{\sum_{\alpha_1, \ldots, \alpha_k=1}^{J} \sum_{i=1}^{J} \vert D_{\alpha (\Gamma )}^{k'}U_i\vert^2 \times \sum_{i=1}^J \vert D_{\alpha (\Gamma^c)}^{k-k'}V_i \vert^2} \\
&=&\vert D^{k'} U\vert \times \vert D^{k-k'}V \vert.
\end{eqnarray*}
This last equality results from the fact that we sum on
different index sets ( $\Gamma $ and  $\Gamma ^{c})$. This gives
\begin{eqnarray*}
\left\vert D^{k} \left\langle U,V\right\rangle_J \right\vert  &\leq
&\sum_{k'=0}^{k}\sum_{\left\vert \Gamma \right\vert =k'}\left\vert
D^{k'}U\right\vert \left\vert D^{k-k'}V\right\vert
=\sum_{k'=0}^{k}C_{k}^{k'}\left\vert D^{k'}U\right\vert \left\vert
D^{k-k'}V\right\vert  \\
&\leq&\sum_{k'=0}^{k}C_{k}^{k'}\left\vert U\right\vert _{k'}\left\vert
V\right\vert _{k-k'} \leq 2^{k}(\sum_{l_{1}+l_{2}=k}\left\vert
U\right\vert _{l_{1}}\left\vert V\right\vert _{l_{2}}).
\end{eqnarray*}
Summing on $k=0,...,l$ we deduce (\ref{scal1}).

\end{pf}
\subsubsection{ Estimation of $\vert \gamma(F)\vert_l$}
We give in this section an estimation of the derivatives of $\gamma(F)$ in terms of $\det \sigma(F)$ and the derivatives of $F$. We assume that $\omega \in \Lambda(F)$.

In what follows $C_{l,d}$ is a constant depending eventually on the order of derivation $l$ and the dimension $d$.
\begin{prop}\label{Bgamma}
Let $F \in \mathcal{S}^d$, we have $\forall l \in \mathbb{N}$

\begin{eqnarray}
\vert \gamma(F) \vert_l & \leq & C_{l,d} \sum_{l_1+l_2 \leq l} \vert F \vert_{1, l_2+1}^{2(d-1)} \left( \frac{1}{\vert \det \sigma(F)\vert} + \sum_{k=1}^{l_1} \frac{\vert F \vert_{1, l_1+1}^{2kd}  }{\vert \det \sigma(F)\vert^{k+1}}\right)   \label{gamma1}\\
& \leq & C_{l,d}  \frac{1}{\vert \det \sigma(F)\vert^{l+1}}(1+\vert F \vert_{1, l+1}^{2d(l+1)}). \label{gamma}
\end{eqnarray}

\end{prop}
Before proving Proposition \ref{Bgamma}, we establish a preliminary lemma.
\begin{lem}\label{lemCR}
 for every $G\in S$, $G>0$ we have
\begin{equation}
\left\vert \frac{1}{G}\right\vert _{l}\leq C_{l}\left(\frac{1}{G}+\sum_{k=1}^{l}
\frac{1}{G^{k+1}}
\sum_{\substack{ k\leq r_{1}+...+r_{k}\leq l \\ 
r_{1},...,r_{k}\geq 1}}
\prod_{i=1}^{k}\left\vert D^{r_i} G\right\vert \right)  
\leq C_{l}(\frac{1}{G}+\sum_{k=1}^{l}\frac{1}{G^{k+1}}\left\vert
G\right\vert _{1,l}^{k}).  \label{boundCR}
\end{equation}
\end{lem}
\begin{pf}
For $F \in \mathcal{S}^d$ and $\phi: \mathbb{R}^d \rightarrow \mathbb{R}$ a $\mathcal{C}^{\infty} $ function, we have from the chain rule
\begin{equation}
D^k_{(\alpha_1, \ldots, \alpha_k)} \phi(F) =\sum_{\vert \beta \vert=1}^k \partial_{\beta} \phi(F) \sum_{\Gamma_1\cup \ldots \cup\Gamma_{\vert \beta \vert}=\{1, \ldots, k\}}\left(
\prod_{i=1}^{ \vert \beta \vert}  D_{\alpha(\Gamma_i)}^{\vert \Gamma_i \vert } F^{\beta_{i}} \right), \label{CRk}
\end{equation}
where $\beta \in \{1, \ldots, d \}^{\vert \beta \vert}$ and $\sum_{\Gamma_1\cup \ldots \cup\Gamma_{\vert \beta \vert}}$ denotes the sum over 
all partitions of $\{1, \ldots,k\}$ with length $\vert \beta \vert$. In particular, for $G \in \mathcal{S}$, $G > 0$ and for $\phi(x)=1/x$, we obtain
\begin{equation}
\vert D^k_{\alpha}( \frac{1}{G}) \vert \leq C_k \sum_{k'=1}^k \frac{1}{ G^{k'+1}} \sum_{\Gamma_1\cup \ldots \cup\Gamma_{k'}=\{1, \ldots, k\}}\left(
\prod_{i=1}^{ k'} \vert D_{\alpha(\Gamma_i)}^{\vert \Gamma_i \vert } G \vert \right). \label{BCRD}
\end{equation}
We deduce then that
\begin{eqnarray*}
\vert D^k(\frac{1}{G} )\vert 
& \leq &
 C_k \sum_{k'=1}^k \frac{1}{ G^{k'+1}} \sum_{\Gamma_1\cup \ldots \cup\Gamma_{k'}=\{1, \ldots, k\}}\left\vert
\prod_{i=1}^{ k'}  D_{\alpha(\Gamma_i)}^{\vert \Gamma_i \vert } G  \right\vert_{\mathbb{R}^{J \otimes k}} , \\
&=& C_k \sum_{k'=1}^k \frac{1}{ G^{k'+1}} \sum_{\Gamma_1\cup \ldots \cup\Gamma_{k'}=\{1, \ldots, k\}}\left(
\prod_{i=1}^{ k'}  \vert D^{\vert \Gamma_i \vert } G\vert \right) , \\
&=& C_k \sum_{k'=1}^k \frac{1}{ G^{k'+1}} 
\sum_{\substack{  r_{1}+...+r_{k'}=k \\ 
r_{1},...,r_{k'}\geq 1}}
\left(\prod_{i=1}^{ k'}  \vert D^{r_i} G\vert \right) , 
\end{eqnarray*}
and the first part of (\ref{boundCR}) is proved. The proof of the second part is straightforward.
\end{pf}

With this lemma, we can prove Proposition \ref{Bgamma}.
\begin{pf}{\bf Proposition \ref{Bgamma}.}
We have on $\Lambda(F)$ 
$$
\gamma^{r,r'}(F)= \frac{1}{\det \sigma(F)} \hat{\sigma}^{r,r'}(F),
$$
where $\hat{\sigma}(F)$ is the algebraic complement of $\sigma(F)$.
But recalling that $\sigma^{r,r'}(F)= \left\langle D^rF,D^{r'} F \right\rangle_J$ we have
\begin{equation}
\vert \det \sigma(F) \vert_l \leq C_{l,d} \vert F \vert_{1,l+1}^{2d} \quad \mbox{ and } \quad
\vert \hat{\sigma}(F) \vert_l \leq C_{l,d} \vert F \vert_{1,l+1}^{2(d-1)}. \label{det}
\end{equation}
Applying inequality (\ref{prod}), this gives
$$
\vert \gamma(F)\vert_l \leq C_{l,d} \sum_{l_1+l_2 \leq l} \vert (\det \sigma(F))^{-1} \vert_{l_1}\vert \hat{\sigma}(F) \vert_{l_2}. 
$$
From Lemma \ref{lemCR} and (\ref{det}), we have
$$
\vert (\det \sigma(F))^{-1} \vert_{l_1} \leq C_{l_1} \left( \frac{1}{\vert \det \sigma(F)\vert} + \sum_{k=1}^{l_1}
\frac{ \vert F \vert_{1, l_1+1}^{2kd}}{ \vert \det \sigma(F) \vert^{k+1}} \right).
$$
Putting together these inequalities, we obtain the inequality (\ref{gamma1}) and consequently (\ref{gamma}).
\end{pf}

\subsubsection{ Some bounds on $H^q$}
Now our goal is to establish some estimates for the weights $H^q$ in terms of the
derivatives of $G$, $F$, $LF$ and $\gamma(F)$. 

\begin{theo} \label{thHq}
For $F \in \mathcal{S}^d$ , $G \in \mathcal{S}$ and  for all
 $q\in \mathbb{N}^*$ there exists an universal constant $C_{q,d}$ such that
for every multi-index $\beta =(\beta _{1},..,\beta _{q})$ 
\begin{eqnarray*}
\left\vert H_{\beta }^{q}(F,G)\right\vert & \leq & \frac{C_{q,d}\left\vert
G\right\vert _{q}(1+\left\vert F\right\vert _{q+1})^{(6d+1)q}}{\left\vert \det
\sigma (F)\right\vert ^{3q-1}}(1+\sum_{j=1}^{q}\sum_{k_{1}+...+k_{j}\leq
q-j}\prod_{i=1}^{j}\left\vert L(F)\right\vert _{k_{i}}), \\
& \leq &\frac{C_{q,d}\left\vert
G\right\vert _{q}(1+\left\vert F\right\vert _{q+1})^{(6d+1)q}}{\left\vert \det
\sigma (F)\right\vert ^{3q-1}}(1+ \vert LF \vert_{q-1}^q).
\end{eqnarray*}
\end{theo}

\begin{pf} 
For $F \in \mathcal{S}^d$, we define the linear operator $T_{r}:
\mathcal{S}\rightarrow \mathcal{S},r=1,...,d$ by 
$$
T_{r}(G)=\left\langle
DG,(\gamma (F)DF)^{r}\right\rangle,
$$
 where $(\gamma (F)DF)^{r} = \sum_{r'=1}^d \gamma^{r',r}(F) DF^{r'}$. Notice that 
\begin{equation}
T_{r}(G\times G^{\prime })=GT_{r}(G^{\prime })+G^{\prime }T_{r}(G).
\label{T1}
\end{equation}
Moreover, for a multi-index $\beta =(\beta _{1},..,\beta _{q})$ we define by
induction $T_{\beta }(G)=T_{\beta _{q}}(T_{(\beta _{1},...,\beta
_{q-1})}(G)).$ We also make the convention that if $\beta $ is the void
multi-index, then $T_{\beta }(G)=G.$ Finally we denote by $L_{r}^{\gamma
}(F)=\sum_{r'=1}^d \delta (\gamma^{r',r} (F)DF^{r'}).$ With this notation we have 
\begin{eqnarray*}
H_{r}(F,G) &=&GL_{r}^{\gamma }(F)-T_{r}(G),  \\
H_{\beta }^{q}(F,G) &=&H_{\beta _{1}}(F,H_{(\beta _{2},...,\beta
_{q})}^{q-1}(F,G)).
\end{eqnarray*}
We will now give an explicite expression of $H_{\beta }^{q}(F,G)$. In
order to do this we have to introduce some more notation. Let $\Lambda_j =\{\lambda_1, \ldots, \lambda_j \} \subset \{1, \ldots, q\}$ such that $\vert \Lambda_j \vert=j$.
We denote by $\mathcal{P}(\Lambda_j)$ the set of the partitions $\Gamma =(\Gamma
_{0},\Gamma _{1},...,\Gamma _{j})$ of $\{1,...,q\} \setminus \Lambda_j$ . Notice that we accept that $\Gamma _{i},i=0,1,...,j$ may be void
sets. Moreover, for a multi-index $\beta =(\beta _{1},..,\beta _{q})$ we
denote by $\Gamma _{i}(\beta )=(\beta _{k_{i}^{1}},...,\beta _{k_{i}^{p}})$
where $\Gamma _{i}=\{k_{i}^{1},...,k_{i}^{p}\}.$ With this notation we can prove by induction and using (\ref{T1}) that 
\begin{equation}
H_{\beta }^{q}(F,G)=T_{\beta }(G)+\sum_{j=1}^{q}\sum_{\Lambda_j \subset \{1, \ldots q\}} \sum_{\Gamma \in \mathcal{P
}(\Lambda_j)}c_{\beta ,\Gamma }T_{\Gamma _{0}(\beta )}(G)\prod_{i=1}^{j}T_{\Gamma
_{i}(\beta )}(L_{\beta_{\lambda_i}}^{\gamma }(F))  \label{T2}
\end{equation}
where $c_{\beta ,\Gamma }\in \{-1,0,1\}$. 

We first give an estimation of $\vert T_{\beta}(G) \vert_l$, for $l \geq 0$ and $\beta=(\beta_1, \ldots, \beta_q)$.
We proceed by induction. For $q=1$ and $1 \leq r \leq d$, we have
$$
\vert T_r(G) \vert_l= \vert \left\langle
DG,(\gamma (F)DF)^{r}\right\rangle \vert_l
$$
and using (\ref{scal3}) we obtain 
\begin{equation*}
\left\vert T_{r}(G)\right\vert _{l} \leq C_l \sum_{l_1+l_2+l_3\leq l}  \vert \gamma(F)\vert_{l_1} \vert G\vert_{l_2+1}  \vert F\vert_{l_3+1}\leq \vert G\vert_{l+1}  \vert F\vert_{l+1} \sum_{l_1=0}^l \vert \gamma(F) \vert_{l_1},
\end{equation*}
where $C_{l}$ is a constant which depends on $l$ only. We obtain then by induction
for every multi-index $\beta=(\beta_1, \ldots, \beta_q) $
\begin{equation}
\left\vert T_{\beta }(G)\right\vert _{l}\leq C_{l,q} \vert G\vert_{l+q}  \vert F\vert^q_{l+q} \sum_{l_1+ \ldots+ l_q\leq l+q-1} \prod_{i=1}^q \vert \gamma(F) \vert_{l_i}. \label{Tbeta}
\end{equation}
In particular this gives for $l=0$ 
$$
\left\vert T_{\beta }(G)\right\vert \leq C_{q} \vert G\vert_{q}  \vert F\vert^q_{q} P_q(\gamma(F)),
$$
with
$$
P_q(\gamma(F))=\sum_{l_1+ \ldots+ l_q\leq q-1} \prod_{i=1}^q \vert \gamma(F) \vert_{l_i}, \quad q\geq1.
$$
To complete the notation, we note $P_0(\gamma(F))=1$. We obtain
\begin{eqnarray*}
\left\vert T_{\Gamma _{i}(\beta )}(L_{\beta_{\lambda_i}}^{\gamma
}(F))\right\vert  
&\leq &C_{q}\left\vert L^{\gamma
}_{\beta_{\lambda_i}}(F)\right\vert _{\left\vert \Gamma _{i}(\beta )\right\vert }
 \vert F\vert^{\left\vert \Gamma _{i}(\beta )\right\vert }_{\left\vert \Gamma _{i}(\beta )\right\vert } P_{\left\vert \Gamma _{i}(\beta )\right\vert }(\gamma(F)).
\end{eqnarray*}

We turn now to the estimation of $\vert L^{\gamma
}_{r}(F)\vert_l$. From the properties of the divergence operator $\delta$ (see Lemma \ref{CR})
$$
\delta (\gamma (F)DF) =\gamma
(F)\delta (DF)-\left\langle D\gamma (F),DF\right\rangle_J .
$$
It follows from (\ref{prod}) and (\ref{scal1}) that
\begin{eqnarray*}
\vert L^{\gamma
}_{r}(F)\vert_l 
\leq C_l \left\vert \gamma (F)\right\vert
_{l+1}(\left\vert \delta (DF)\right\vert _{l}+\left\vert
F\right\vert _{l+1}) \leq 
C_l \left\vert \gamma (F)\right\vert
_{l+1}(1+\left\vert LF\right\vert _{l})(1+\left\vert
F\right\vert _{l+1}),
\end{eqnarray*}
and we get 
\begin{eqnarray}
\left\vert T_{\Gamma _{i}(\beta )}(L_{\beta_{\lambda_i}}^{\gamma
}(F))\right\vert  
\leq C_{q}
\left\vert \gamma (F)\right\vert
_{\left\vert \Gamma _{i}(\beta )\right\vert+1}(1+\left\vert LF\right\vert _{\left\vert \Gamma _{i}(\beta )\right\vert})(1+\left\vert
F\right\vert _{\left\vert \Gamma _{i}(\beta )\right\vert+1})
 \vert F\vert^{\left\vert \Gamma _{i}(\beta )\right\vert }_{\left\vert \Gamma _{i}(\beta )\right\vert } P_{\left\vert \Gamma _{i}(\beta )\right\vert }(\gamma(F)).
\end{eqnarray}
Reporting these inequalities in (\ref{T2}) and  recalling that $\left\vert \Gamma _{0}(\beta )\right\vert + \ldots +
\left\vert \Gamma _{j}(\beta )\right\vert=q-j$ we deduce :
\begin{eqnarray}
\vert H_{\beta}^q(F,G)\vert& \leq & \vert T_{\beta}(G) \vert + 
C_{q,d} \sum_{j=1}^q \sum_{k_0+ \dots + k_j=q-j} 
\vert G \vert_{k_0} \vert F \vert_{k_0}^{k_0} P_{k_0}(\gamma(F)) \left( \prod_{i=1}^j \vert \gamma(F) \vert_{k_i+1} P_{k_i}(\gamma(F)) \right.  \nonumber \\
 & & \quad \quad  \quad \quad \quad \left. \vert F \vert_{k_i}^{k_i}(1 + \vert F \vert_{k_i+1})(1+\vert LF \vert_{k_i})\right) \label{T3}
\end{eqnarray}

Now, for $q \geq 1$, we have from (\ref{gamma}) :
$$
P_q(\gamma(F)) \leq C_q \frac{1}{ \vert \det \sigma(F) \vert^{2q-1}}(1 + \vert F \vert_q)^{4dq}, 
$$
so the following inequality holds for $q\geq 0$ :
$$
P_q(\gamma(F)) \leq C_q \frac{1}{ \vert \det \sigma(F) \vert^{2q}}(1 + \vert F \vert_q)^{4dq}.
$$
We obtain then for
$k_0,k_1, \ldots, k_j \in \mathbb{N}$ such that $k_0+ \ldots + k_j=q-j$ 
\begin{eqnarray}
\prod_{i=0}^j P_{k_i}(\gamma(F)) \leq C_q \frac{1}{ \vert \det \sigma(F) \vert^{2(q-j)}}(1 + \vert F \vert_{q-j})^{4d(q-j)} 
\end{eqnarray}
and once again from (\ref{gamma})
\begin{eqnarray}
\prod_{i=1}^j \vert \gamma(F) \vert_{k_i+1} \leq C_q \frac{1}{ \vert \det \sigma(F) \vert^{q+j}}(1 + \vert F \vert_{q-j+2})^{2d(q+j)}
\end{eqnarray}
it yields finally 
$$
\prod_{i=0}^j P_{k_i}(\gamma(F)) \prod_{i=1}^j \vert \gamma(F) \vert_{k_i+1} \leq C_q \frac{1}{ \vert \det \sigma(F) \vert^{3q-j}}(1 + \vert F \vert_{q-j+2})^{6dq-2dj}.
$$
Turning back to  (\ref{T3}),  it
follows that 
\begin{eqnarray*}
\left\vert H_{\beta }^{q}(F,G)\right\vert & \leq & \frac{C_{q,d}\left\vert
G\right\vert _{q}(1+\left\vert F\right\vert _{q+1})^{(6d+1)q}}{\left\vert \det
\sigma (F)\right\vert ^{3q-1}}(1+\sum_{j=1}^{q}\sum_{k_{1}+...+k_{j}\leq
q-j}\prod_{i=1}^{j}\left\vert L(F)\right\vert _{k_{i}}),
\end{eqnarray*}
and Theorem \ref{thHq} is proved.
\end{pf}

\section{Stochastic equations with jumps}

\subsection{Notations and hypotheses}

We consider a Poisson point process $p$ with state space $(E,B(E)),$ where $
E=\mathbb{R}^{d}\times \mathbb{ R}_{+}.$\ We refer to [I.W] for the notation. We denote by $N$
the counting measure associated to $p$, we have  $N([0,t)\times A)=\# \{ 0 \leq s <t; p_s \in A \}$ for $t \geq 0$ and $A \in B(E)$. 
We assume that the associated intensity measure is given by $\widehat{
N}(dt,dz,du)=dt\times d\mu (z)\times 1_{[0,\infty )}(u)du$ where $(z,u)\in
E=\mathbb{R}^{d}\times \mathbb{R}_{+}$ and $\mu (dz)=h(z)dz.$

We are interested  in the solution of the $d$ dimensional stochastic equation
\begin{equation}
X_{t}=x+\int_{0}^{t}\int_{E}c(z,X_{s-})1_{\{u<\gamma
(z,X_{s-})\}}N(ds,dz,du)+\int_{0}^{t}g(X_{s})ds.  \label{eq1}
\end{equation}
We remark that the infinitesimal generator of
the Markov process $X_{t}$ is given by 
\begin{equation*}
L\psi (x)=g(x)\nabla \psi (x)+\int_{\mathbb{R}^d}(\psi (x+c(z,x))-\psi (x))K(x,dz)
\end{equation*}
where $K(x,dz)=\gamma (z,x)h(z)dz$  depends on the variable $x\in
\mathbb{R}^{d}. $\ See [F.1]  for the proof of existence and
uniqueness of the solution of the above equation.

Our aim is to give sufficient conditions in order to prove that the law of $X_{t}$ is
absolutely continuous with respect to the Lebesgue measure and has a smooth
density. In this section we make the following hypotheses on the
functions $\gamma, g, h$ and $c$.

{\bf Hypothesis 3.0}  We assume that $\gamma ,g,h$
and $c$ are infinitely differentiable functions in both variables $z$ and $
x$. Moreover we assume that $g$ and its derivatives are bounded and that $\ln h$ has bounded derivatives

\textbf{Hypothesis 3.1}. We assume that there exist two functions $
\overline{\gamma },\underline{\gamma }:\mathbb{R}^{d}\rightarrow \mathbb{R}_{+}$ such that
\begin{equation*}
\overline{C}\geq \overline{\gamma }(z)\geq \gamma (z,x)\geq 
\underline{\gamma }(z)\geq 0,\quad \forall x\in \mathbb{R}^{d}
\end{equation*}
where $\overline{C}$ is a constant. 

\textbf{Hypothesis 3.2}. \textbf{i)} We assume that there exists a non negative
and bounded function $\overline{c}:\mathbb{R}^{d}\rightarrow \mathbb{R}_{+}$ such that $
\int_{\mathbb{R}^{d}}\overline{c}(z)d\mu (z)<\infty $ and 
\begin{equation*}
\left\vert c(z,x)\right\vert +\left\vert \partial _{z}^{\beta }\partial
_{x}^{\alpha }c(z,x)\right\vert \leq \overline{c}(z)\quad \forall z,x\in
\mathbb{R}^{d}.
\end{equation*}
We need this hypothesis  in order to estimate the Sobolev norms.

\textbf{ii)} There exists a measurable function $\widehat{c}
:\mathbb{R}^{d}\rightarrow \mathbb{ R}_{+}$ such that $\int_{\mathbb{R}^{d}}\widehat{c}(z)d\mu (z)<\infty $ and
\begin{equation*}
\left\Vert \nabla _{x}c\times (I+\nabla _{x}c)^{-1}(z,x)\right\Vert \leq 
\widehat{c}(z),\quad \forall (z,x)\in \mathbb{R}^{d}\times \mathbb{R}^{d}.  
\end{equation*}
In order to simplify the notations we assume that $\widehat{c}(z)=\overline{c}(z).$

\textbf{iii) }There exists a non negative function $\underline{c}
:\mathbb{R}^{d}\rightarrow \mathbb{R}_{+}$ such that for every $z\in \mathbb{R}^{d}$
\begin{equation*}
\sum_{r=1}^{d}\left\langle \partial _{z_{r}}c(z,x),\xi \right\rangle
^{2}\geq \underline{c}^{2}(z)\left\vert \xi \right\vert ^{2},\quad \forall
\xi \in \mathbb{R}^{d}  
\end{equation*}
and we assume that there exists $\theta >0$ such that 
\begin{equation*}
\underline{\lim }_{a\rightarrow + \infty} \frac{1}{\ln a} \int_{\{\underline{c}^{2}\geq1/ a\}}
\underline{\gamma }(z)d\mu (z)= \theta.
\end{equation*}
{\bf Remark} : assumptions ii) and iii) give sufficient conditions   to prove the non degeneracy of the Malliavin covariance matrix as defined in the previous section.
In particular the second part of iii) implies that $\underline{c}^2$ is a $(p,t)$ broad function (see $[B.G.J.]$) for $p/t< \theta$.
Notice that we may have $\underline{c}(z)=0$ for some $z\in \mathbb{R}^{d}.$

We add to these hypotheses some assumptions on the derivatives of $\gamma$ and $\ln \gamma$ with respect to $x$ and $z$. For $l \geq 1$ we use the notation :
\begin{eqnarray*}
\overline{\gamma}^{x,l}(z)=\sup_x \sup_{1 \leq \vert \beta \vert \leq l} \vert \partial_{\beta,x} \gamma(z,x) \vert, \\
\overline{\gamma}_{\ln}^{x,l}(z)=\sup_x \sup_{1 \leq \vert \beta \vert \leq l} \vert \partial_{\beta,x} \ln \gamma(z,x) \vert,\\
\overline{\gamma}_{\ln}^{z,l}(z)=\sup_x \sup_{1 \leq \vert \beta \vert \leq l} \vert \partial_{\beta,z} \ln \gamma(z,x) \vert.
\end{eqnarray*}
\textbf{Hypothesis 3.3}. We assume that $\ln \gamma$ has bounded derivatives with respect to $z$  (that is $\overline{\gamma}_{\ln}^{z,l}(z)$ is bounded) and that  $\gamma$ has bounded derivatives with respect to $x$ such that $\forall z \in \mathbb{R}^d $, $\overline{\gamma}^{x,l}(z) \leq \overline{\gamma}^{x,l};$ moreover we assume that
$$
\sup_{z^* \in \mathbb{R}^d} \int_{B(z^*,1)} \overline{\gamma}(z) d \mu(z)<+ \infty.
$$
We complete this hypothesis with two alternative hypotheses.

a) (weak dependence on $x$) We assume that $\forall l \geq 1$
$$
\int_{\mathbb{R}^d} \overline{\gamma}_{\ln}^{x,l}(z) \overline{\gamma}(z) d\mu(z) < \infty .  
$$

b) (strong dependence on $x$) We  assume that  $\ln \gamma$ has bounded derivatives
with respect to $x$ such that $\forall l \geq 1$
$$
\forall z \in \mathbb{R}^d, \quad \overline{\gamma}_{\ln}^{x,l}(z) \leq \overline{\gamma}_{\ln}^{x,l} .
$$
{\bf Remark} : if $\mu$ is the Lebesgue measure ( case $h=1$) and if $\gamma$ does not depend on $z$ then $\overline{\gamma}_{\ln}^{x,l}$ is constant and consequently hypothesis $3.3.a$ fails. Conversely, if $\gamma(z,x)=\gamma(z)$ then hypothesis $3.3.a$ is satisfied as soon as 
$\ln \gamma$ has bounded derivatives. This last case corresponds to the standard case where the law of the amplitude of the jumps does not depend on the position of $X_t$. Under Hypothesis $3.3.a$ we are in  a classical situation where the divergence does not blow up and this leads to an integration by part formula with bounded weights (see Proposition \ref{pIPPFM} and Lemma \ref{LZ}). On the contrary under assumption $3.3.b$, the divergence can blow up as well as the weights appearing in the integration by part formula.

\subsection{Main results and examples}

Our methodology to study the regularity of the law of the random variable $X_t$ is based on the following result.  Let $\hat{p}_X(\xi)=
E(e^{i \left\langle \xi,X \right\rangle})$ be the Fourier transform of a $d$-dimensional random variable $X$ then using the Fourier inversion formula, one can prove
that if  $\int_{\mathbb{R}^d} \vert \xi \vert^p \vert \hat{p}_X(\xi)\vert d \xi < \infty$  for $p>0$ then the law of $X$ is absolutely continuous with respect to the 
Lebesgue measure on $\mathbb{R}^d$ and its density is $\mathcal{C}^{[p]}$, where $[p]$ denotes the entire part of $p$.

To apply this result, we just have to bound the Fourier transform of $X_t$ in terms of $1/\vert \xi \vert.$ This is done in the next proposition. 
The proof of this proposition
needs a lot of steps that we detail in the next sections and it will be given later.

\begin{prop} \label{fourierX}
Let $B_M=\{ z \in \mathbb{R}^d ; \vert z \vert < M \}$, then under hypotheses $3.0.$, $3.1.$ $3.2.$ and $3.3$ we have for all  $M \geq 1$,  for $q \geq 1$ and $t>0$ such that $4d(3q-1)/t< \theta$

a) if $3.3.a$ holds
\begin{eqnarray*}
\vert \hat{p}_{X_t}(\xi) \vert   &\leq & t \int_{B_{M-1}^{c}}\underline{c}^2(z)\underline{\gamma }(z)d\mu (z)  \frac{1}{2}\left\vert \xi \right\vert ^{2}+\left\vert \xi \right\vert t e^{Ct}
\int_{B_{M}^{c}}\overline{c}(z)\overline{\gamma }(z)d\mu (z) 
+\frac{C_q}{\left\vert \xi \right\vert ^{q}}.
\end{eqnarray*}

b) if $3.3.b$ holds

\begin{eqnarray*}
\vert \hat{p}_{X_t}(\xi) \vert   &\leq &t \int_{B_{M-1}^{c}}\underline{c}^2(z)\underline{\gamma }(z)d\mu (z)  \frac{1}{2}\left\vert \xi \right\vert ^{2}+\left\vert \xi \right\vert t e^{Ct}
\int_{B_{M}^{c}}\overline{c}(z)\overline{\gamma }(z)d\mu (z) 
+\frac{C_q(1+\mu(B_{M+1})^q)}{\left\vert \xi \right\vert ^{q}}.
\end{eqnarray*}

\end{prop}
We can remark that if $\theta = +\infty$ then the result holds $\forall q \geq 1$ and $\forall t>0$.

By choosing $M$ judiciously as a function of $\xi$ in
the inequalities given in Proposition \ref{fourierX}, we obtain $ \vert \hat{p}_{X_t}(\xi) \vert \leq C/ \vert \xi \vert^p$ for some $p>0$ and this permits us
 to deduce some regularity for the density of $X_t$. The next theorem precise the optimal choice of $M$
 with respect to $\xi$ and permits us to derive
the regularity of the law of the process $X_t$. 
\begin{theo} \label{density}
We assume that hypotheses $3.0.$, $3.1.$, $3.2$ and $3.3.$ hold.

{\bf a)} Assuming $3.3.a$, the law of $X_t$ admits a density $\mathcal{C}^k$ if  $t>(3k+3d-1) \frac{4d}{\theta}$. In the case $\theta= \infty$, the law of $X_t$ admits a density $\mathcal{C}^{\infty}$.

{\bf b)} Assuming  $3.3.b$ and the two following hypotheses

{\bf A1} : $\exists p_1,p_2>0$ such that :
$$
\limsup_M M^{p_1} \int_{B_M^c} \overline{c}(z) \overline{\gamma}(z) d \mu(z) < + \infty ;
$$
$$
\limsup_M M^{p_2} \int_{B_M^c} \underline{c}^2(z) \underline{\gamma}(z) d \mu(z) < + \infty ;
$$

{\bf A2} : $\exists \rho >0$ such that $\mu(B_M) \leq C M^{\rho}$ where $B_M=\{ z \in \mathbb{R}^d ; \vert z \vert < M \}$;

{\bf case 1}: if  $\theta=+ \infty$ then the law of $X_t$ admits a density $\mathcal{C}^k$ with $k<\min(p_1/\rho-1-d,p_2/\rho-2-d)$ if $\min(p_1/\rho-1-d,p_2/\rho-2-d) \geq 1$ .

{\bf case 2} : if $0< \theta<\infty$ let $q^*(t ,\theta)=[ \frac{1}{3}(\frac{t \theta}{4d}+1) ]$;  then
the law of $X_t$ admits a density $\mathcal{C}^k$ for $k<\sup_{0<r<1/\rho}\min(r p_1-1-d, r p_2-2-d,q^*(t, \theta)(1-r \rho)-d)$,
 if for some $0<r<1/\rho$, $\min(r p_1-1-d, r p_2-2-d,q^*(t, \theta)(1-r \rho)-d)\geq 1$.
\end{theo}

\begin{pf} 

a) Assuming $3.3.a$
and letting $M$ go to infinity in the right-hand side of the inequality given in Proposition \ref{fourierX} , we deduce
$$
\vert \hat{p}_{X_t}(\xi) \vert  \leq C/ \vert \xi \vert^q,
$$ 
and the result follows.

b) From $A1$, for $M$ large enough, we have
$$
\int_{B_{M}^{c}}\overline{c}(z)\overline{\gamma }(z)d\mu (z) \leq C/M^{p_1}
$$
and
$$
\int_{B_{M-1}^{c}}\underline{c}^2(z)\underline{\gamma }(z)d\mu (z) \leq C/M^{p_2}.
$$
Now assuming $3.3.b$ and $A2$ and 
choosing $M= \vert \xi \vert^r$, for $0<r<1/ \rho$, we obtain from Proposition \ref{fourierX} 
\begin{eqnarray*}
\vert \hat{p}_{X_t}(\xi) \vert \leq C \left(\frac{1}{ \vert \xi \vert^{rp_1-1}}+ \frac{1}{ \vert \xi \vert^{rp_2-2}}+ \frac{1}{\left\vert \xi \right\vert ^{q(1-r\rho)}}\right),
\end{eqnarray*}
for $q$ and $t$ such that $4d(3q-1)/t< \theta$. Now if $\theta= \infty$, we obtain for $q$ large enough
$$
\vert \hat{p}_{X_t}(\xi) \vert \leq C \left(\frac{1}{ \vert \xi \vert^{rp_1-1}}+ \frac{1}{ \vert \xi \vert^{rp_2-2}} \right).
$$
In the case $\theta<\infty$, the best choice of $q$ is $q^*(t ,\theta)$.
This achieves the proof of theorem \ref{density}.
\end{pf}

We end this section with some examples in order to illustrate the results of Theorem \ref{density}. 

\textbf{Example 1}. In this example we assume that $h=1$ so $\mu(dz)=dz$ and that $\underline{\gamma }(z)$
is equal to a constant $\underline{\gamma}>0$. We also assume that Hypothesis $3.3 .b$ holds.
We have $\mu (B_{M})=r_{d}M^{d}$
where $r_{d}$ is the volume of the unit ball in $\mathbb{R}^{d}$ so $\rho =d.$\ We
will consider two types of behaviour for $c.$

\textbf{i)} \textbf{Exponential decay}: we assume that $\overline{c}
(z)=e^{-b\left\vert z\right\vert ^{c}}$ and $\underline{c}
(z)=e^{-a\left\vert z\right\vert ^{c}}$ for some constants $0<b\leq a$ and $
c>0.$ We have
\begin{equation*}
\int_{\{\underline{c}^{2}>1/u\}}\underline{\gamma }(z)d\mu (z)=\frac{
\underline{\gamma}r_{d}}{(2a)^{d/c}}\times (\ln u)^{d/c}.
\end{equation*}
We deduce then
\begin{equation}
\theta =0\quad if\quad c>d,\quad \theta =\infty \quad if\quad 0<c<d\quad
and\quad \theta =\frac{\underline{\gamma}r_{d}}{2a}\quad if\quad c=d.  \label{ex1}
\end{equation}
If $c>d$, hypothesis $3.2.iii$ fails,  this is coherent with the result
of $[B.G.J].$ Now observe that 
\begin{equation*}
\int_{B_{M}^{c}}\underline{c}^{2}(z)\underline{\gamma }(z)d\mu
(z)+\int_{B_{M}^{c}}\overline{c}(z)\overline{\gamma }(z)d\mu (z)\leq
e^{-\eta \left\vert z\right\vert ^{c}}
\end{equation*}
for some $\eta >0$ so $p_{1}=p_{2}=\infty .$ In the case $0<c<d$ we obtain a
density  $C^{\infty }$ for every $t>0.$ In the case $c=d$ we have $
q^{\ast }(t,\theta )=[\frac{1}{3}(1+\frac{\underline{\gamma}r_{d}}{8da}\times t)].$
If $t < 8da(3d+2)/(\underline{\gamma}r_{d})$ we obtain nothing and if $
t \geq 8da(3d+2)/(\underline{\gamma}r_{d})$ we obtain a density  $\mathcal{C}^k$ where $k$
is the largest integer less than $[\frac{1}{3}(1+\frac{\underline{\gamma}r_{d}}{
8da}\times t)]-d.$

\textbf{ii)} \textbf{Polynomial decay}. We assume that $\overline{c}
(z)=b/(1+\left\vert z\right\vert ^{p})$ and $\underline{c}
(z)=a/(1+\left\vert z\right\vert ^{p})$ for some constants $0<a\leq b$ and $
p>d.$ We have
\begin{equation*}
\int_{\{\underline{c}^{2}>1/u\}}\underline{\gamma }(z)d\mu (z)=\underline{\gamma}
r_{d}\times (a\sqrt{u}-1)^{d/p}
\end{equation*}
so  $\theta =\infty $ and our result works for every $t>0.$ Hence a simple computation gives
\begin{equation*}
\int_{B_{M}^{c}}\underline{c}^{2}(z)\underline{\gamma }(z)d\mu (z)\leq 
\frac{C}{M^{2p-d}},\quad \int_{B_{M}^{c}}\overline{c}(z)\overline{\gamma }
(z)d\mu (z)\leq \frac{C}{M^{p-d}}
\end{equation*}
and then $p_{1}=p-d$ and $p_{2}=2p-d.$ If  
$p\geq d(d+3)$ then $\min( p/d-2-d,2p/d-3-d) \geq 1$ and we obtain a density $\mathcal{C}^k$ with $k<\frac{p}{d}-d-2.$
Conversely if $p<d(d+3)$, we can say nothing about the regularity of the density of $X_t$.
We give now an example where the function $\gamma$ satisfies Hypothesis $3.3. a $. 

\textbf{Example 2}. As in the preceding example, we assume $h=1$. We consider the function $\gamma (z,x)=\exp (-\alpha
(x)/(1+\left\vert z\right\vert ^{q}))$ for some $q>d.$ We assume that $
\alpha $ is a smooth function which is bounded and has bounded derivatives
and moreover there exists two constants such that $\overline{\alpha }\geq
\alpha (x)\geq \underline{\alpha }>0.$ Notice that the derivatives with
respect to $x$ of $\ln \gamma (z,x)$ are bounded by $C/(1+\left\vert
z\right\vert ^{q})$ which is integrable with respect to the Lebesgue measure
if $q>d.$ So Hypothesis $3.3. a$ is true. Moreover we check that $\underline{\gamma}(z)=\exp (-\overline{\alpha}
/(1+\left\vert z\right\vert ^{q}))$.

\textbf{i)} \textbf{Exponential decay}. We take $c$ as in
Example 1.i). It follows that 
\begin{equation*}
 \int_{\{
\underline{c}^{2}>1/u\}}\underline{\gamma }(z)d\mu (z)\geq \exp (-\overline{\alpha}) \frac{r_{d}}{(2a)^{d/c}}\times (\ln u)^{d/c}.
\end{equation*}
So we obtain once again $\theta $ as in (\ref{ex1}). In the
case $c>d$ we can say nothing, in the case $c<d$ we obtain a density 
$C^{\infty }$ and in the case $c=d$ we have $\theta= \frac{r_{d}}{2a}$ and we obtain a density  $C^{k}$ if $
t>\frac{8ad(3k+3d-1)}{r_{d}}.$ In particular we have no results
if $t\leq \frac{8ad(3d-1)}{r_{d}}.$ Notice that the only
difference with respect to the previous example concerns the case $c=d$ when
we have a slight gain.

\textbf{ii)} \textbf{Polynomial decay}. At last we take $c
$ as in the example 1.ii). We check that $\theta =\infty $ so we obtain a density  $C^{\infty }$
, which is a better result than the one of the previous example.

{\bf Example 3.}
We consider the process $(Y_t)$ solution of the stochastic equation 
\begin{equation*}
dY_{t}=f(Y_{t})dL_{t},
\end{equation*}
where $L_t$ is a L\'{e}vy process with intensity measure $\vert y \vert^{-(1+\rho
)}1_{\{\left\vert y\right\vert \leq 1\}}dy$, with $0<\rho<1$. 
The infinitesimal generator of $Y$ is  given by
\begin{equation*}
L\psi (x)=\int_{\{\left\vert y\right\vert \leq 1\}}(\psi (x+f(x)y)-\psi (x))
\frac{dy}{\vert y \vert^{1+\rho }}.
\end{equation*}
If we introduce some function $g(x)$ in this operator we obtain 
\begin{equation*}
L\psi (x)=\int_{\{\left\vert y\right\vert \leq 1\}}(\psi (x+f(x)y)-\psi
(x))g(x)\frac{dy}{\vert y \vert^{1+\rho }}.
\end{equation*}
We are interested to represent this operator through a stochastic equation. In order
to come back in our framework, we translate the integrability problem from $0$ to $\infty$ by  the change of variables $z=y^{-1}$
and we obtain 
\begin{equation*}
L\psi (x)=\int_{\{\left\vert z\right\vert \geq 1\}}(\psi (x+f(x)z^{-1})-\psi
(x))g(x)\frac{dz}{\vert z \vert^{1-\rho }}.
\end{equation*}
This operator can be viewed as the infinitesimal generator of the process $(X_t)$ solution of
$$
X_t=x+ \int_0^t \int_{\mathbb{R} \times \mathbb{R}_+} f(X_{s-})z^{-1} 1_{\{u<g(X_{s-})\}} N(ds, dz,du).
$$
We have $E=\mathbb{R} \times \mathbb{R}_+$, $d\mu(z)=\frac{1}{\vert z \vert^{1-\rho }} 1_{\{\vert z\vert \geq 1\}}dz$, $c(z,x)=f(x)z^{-1}$ and
$\gamma (z,x)=g(x).$
We make the following assumptions. There exist two constants $\underline{f}$ and $\overline{f}$ such that $\forall x$ $\underline{f} \leq f(x) \leq \overline{f}$ and  we suppose that all derivatives of $f$ are bounded by $\overline{f}$. Moreover we assume that  there exist two constants $\underline{g}$ and $\overline{g}$  such that $g$ and its derivative are bounded by $\overline{g}$ and $0<\underline{g}  \leq g(x)$, $\forall x$.
Consequently it is easy to check that hypotheses  $3.0.$, $3.1.$, $3.2.$ and $3.3. b$ are satisfied, with $\theta=+ \infty$. Moreover we have $\mu(B_M) \leq C M^{\rho}$ and $A2$ holds  with $p_1=1-\rho$ and $p_2=2 -\rho$.
Consequently we deduce that the law of $X_t$ admits a density $\mathcal{C}^k$ with $k< 1/\rho -3$ if $1/ \rho-3 \geq 1$.

The next sections are the successive steps to prove proposition \ref{fourierX} .

\subsection{Approximation of $X_t$}

In order to prove that the process $X_t$, solution of (\ref{eq1}),  has a smooth density, we will apply the differential calculus and the integration by parts formula of section $2$. But since the random variable $X_t$ can not be viewed as a simple functional, the first step consists in approximate it. We describe in this section our approximation procedure. We consider a non-negative and smooth function $
\varphi :\mathbb{R}^{d}\rightarrow \mathbb{R}_{+}$ such that $\varphi (z)=0$ for $\left\vert
z\right\vert >1$ and $\int_{\mathbb{R}^{d}}\varphi (z)dz=1.$ And for $M\in \mathbb{N}$ we
denote $\Phi _{M}(z)=\varphi \ast 1_{B_{M}}$ with $B_{M}=\{z\in
\mathbb{R}^{d}:\left\vert z\right\vert <M\}.$ Then $\Phi _{M}\in C_{b}^{\infty }$ and
we have $1_{B_{M-1}}\leq \Phi _{M}\leq 1_{B_{M+1}}.$ We denote by $X_{t}^{M}$
the solution of the equation 
\begin{equation}
X_{t}^{M}=x+\int_{0}^{t}\int_{E}c_{M}(z,X_{s-}^{M})1_{\{u<\gamma
(z,X_{s-}^{M})\}}N(ds,dz,du)+\int_{0}^{t}g(X_{s}^{M})ds. \label{eq2}
\end{equation}
where $
c_{M}(z,x):=c(z,x)\Phi _{M}(z)$. Observe that equation $(\ref{eq2})$ is obtained from $(\ref{eq1})$ replacing the coefficient $c$ by the truncating one $c_M$.
Let $N_{M}(ds,dz,du):=1_{B_{M+1}}(z)\times 1_{[0,2 \overline{C}
]}(u)N(ds,dz,du).$ Since $\{u<\gamma (z,X_{s-}^{M})\}\subset \{u<2 \overline{C}
\}$ and $\Phi _{M}(z)=0$ for $\left\vert z\right\vert >M+1,$ we may replace $
N$ by $N_{M}$ in the above equation and consequently $X_t^M$ is solution of the equation
\begin{equation*}
X_{t}^{M}=x+\int_{0}^{t}\int_{E}c_{M}(z,X_{s-}^{M})1_{\{u<\gamma
(z,X_{s-}^{M})\}}N_{M}(ds,dz,du)+\int_{0}^{t}g(X_{s}^{M})ds.
\end{equation*}

Since the intensity measure $\widehat{N}_{M}$ is finite we may
represent the random measure  $N_{M}$ by a compound Poisson process. Let $
\lambda_M = 2 \overline{C}\times \mu (B_{M+1})=t^{-1}E(N_{M}(t,E))$ and let
 $J_{t}^M$ a Poisson process of parameter $\lambda_M .$ We denote by $
T_{k}^M,k\in \mathbb{N}$ the jump times of $J_{t}^M$.  We also consider two sequences of
independent random variables $(Z_{k}^M)_{k \in \mathbb{N}}$ and $(U_{k})_{k \in \mathbb{N}}$ respectively in $\mathbb{R}^d$ and $\mathbb{ R}_{+}$ 
which are independent of $J^M$ and such that
\begin{equation*}
Z_{k}\sim \frac{1}{\mu (B_{M+1})}1_{B_{M+1}}(z)d\mu(z),\quad and\quad U_{k}\sim 
\frac{1}{2 \overline{C}}1_{[0,2 \overline{C}]}(u)du.
\end{equation*}
To simplify the notation, we omit the dependence on $M$ for the variables $(T_k^M)$ and $(Z_k^M)$.
Then  equation $(\ref{eq2})$ may be written as 
\begin{equation}
X_{t}^{M}=x+\sum_{k=1}^{J_{t}^M}c_{M}(Z_{k},X_{T_{k}-}^{M})1_{(U_{k},\infty
)}(\gamma (Z_{k},X_{T_{k}-}^{M}))+\int_{0}^{t}g(X_{s}^{M})ds.  \label{eq3}
\end{equation}

\begin{lem}\label{lapprox}
Assume that hypotheses $3.0.$, $ 3.1.$, $3.2$ and $3.3.$ hold true
then we have
\begin{equation}
E\left\vert X_{t}^{M}-X_{t}\right\vert \leq \varepsilon
_{M}:=te^{Ct}\int_{\{\left\vert z\right\vert >M\}}\overline{c}(z)\overline{
\gamma }(z)d\mu (z),  \label{approx}
\end{equation}
for some constant $C$.
\end{lem}

\begin{pf} We have $E\left\vert X_{t}^{M}-X_{t}\right\vert \leq
I_{M}^{1}+I_{M}^{2}$ with
\begin{eqnarray*}
I_{M}^{1} &=&E\int_{0}^{t}\int_{\mathbb{R}^{d}}\int_{0}^{\overline{C}}\left\vert
c(z,X_{s})1_{\{u<\gamma (z,X_{s})\}}-c_{M}(z,X_{s}^{M})1_{\{u<\gamma
(z,X_{s}^{M})\}}\right\vert dud\mu (z)ds \\
I_{M}^{2} &=&E\int_{0}^{t}\left\vert g(X_{s})-g(X_{s}^{M})\right\vert ds.
\end{eqnarray*}
Since $\left\vert \nabla _{x}c(z,x)\right\vert \leq \overline{c}(z)$ we have 
$I_{M}^{1}\leq I_{M}^{1,1}+I_{M}^{1,2}$ with 
\begin{eqnarray*}
I_{M}^{1,1} &=&E\int_{0}^{t}\int_{\mathbb{R}^{d}}\int_{0}^{\overline{C}}\left\vert
c(z,X_{s})-c_{M}(z,X_{s}^{M})\right\vert 1_{\{u<\overline{\gamma }
(z)\}}dud\mu (z)ds \\
&\leq &t\int_{\mathbb{R}^{d}}\overline{c}(z)\overline{\gamma }(z)(1-\Phi _{M}(z))d\mu
(z)+\int_{\mathbb{R}^{d}}\overline{c}(z)\overline{\gamma }(z)dz\times
E\int_{0}^{t}\left\vert X_{s}-X_{s}^{M}\right\vert ds
\end{eqnarray*}
and, since $\left\vert \nabla _{x}\gamma (z,x)\right\vert \leq \overline{\gamma}^{x,1}$
\begin{eqnarray*}
I_{M}^{1,2} &=&E\int_{0}^{t}\int_{\mathbb{R}^{d}}\int_{0}^{\overline{C}}\overline{c}
(z)\left\vert 1_{\{u<\gamma (z,X_{s})\}}-1_{\{u<\gamma
(z,X_{s}^{M})\}}\right\vert dud\mu (z)ds \\
&=&E\int_{0}^{t}\int_{\mathbb{R}^{d}}\overline{c}(z)\left\vert \gamma
(z,X_{s})-\gamma (z,X_{s}^{M})\right\vert d\mu (z)ds \\
&\leq &  \int_{\mathbb{R}^{d}}\overline{c}(z) \overline{\gamma}^{x,1} d\mu (z)\times
E\int_{0}^{t}\left\vert X_{s}-X_{s}^{M}\right\vert ds.
\end{eqnarray*}
A similar inequality holds for $I_{M}^{2}$ so we obtain
\begin{equation*}
E\left\vert X_{t}^{M}-X_{t}\right\vert \leq t\times \int_{R^{d}}\overline{
\gamma }(z)\overline{c}(z)(1-\Phi _{M}(z))d\mu (z)+C\int_{0}^{t}E\left\vert
X_{s}-X_{s}^{M}\right\vert ds.
\end{equation*}
We conclude by using Gronwall's lemma. 
\end{pf}

The random variable $X^M_t$ solution of $(\ref{eq3})$ is a function of $(Z_1 \ldots, Z_{J_t^M})$ but it is not a simple functional, as defined
in section $2$
because the coefficient 
$c_{M}(z,x)1_{(u,\infty )}(\gamma (z,x))$ is not differentiable with respect
to $z$. In order to avoid this difficulty we use the following
alternative representation. Let $z_{M}^*\in \mathbb{R}^{d}$ such that $
\left\vert z_{M}^*\right\vert =M+3$.  We define 
\begin{eqnarray}
q_{M}(z,x) &:&=\varphi (z-z_{M}^*)\theta _{M,\gamma }(x)+\frac{1}{2\overline{C}
\mu (B_{M+1})}1_{B_{M+1}}(z)\gamma (z,x)h(z)  \label{qM} \\
\theta _{M,\gamma }(x) &:&=\frac{1}{\mu (B_{M+1})}\int_{\{\left\vert
z\right\vert \leq M+1\}}(1-\frac{1}{2 \overline{C}}\gamma (z,x))\mu (dz). 
\notag
\end{eqnarray}
We recall that $\varphi $\ is the function defined at the beginning of this subsection : a non-negative and smooth
function with $\int \varphi =1$ and which is null outside the unit ball.  Moreover from hypothesis $3.1$, $0 \leq \gamma(z,x) \leq \overline{C}$ and then $1 \geq \theta_{M, \gamma}(x) \geq 1/2$.
By construction
the function $q_M$ satisfies $\int
q_{M}(x,z)dz=1. $ Hence we can check that
\begin{equation}
E(f(X_{T_{k}}^{M})\mid
X_{T_{k}-}^{M}=x)=\int_{R^{d}}f(x+c_{M}(z,x))q_{M}(z,x)dz.  \label{condXM}
\end{equation}
In fact the  left hand side term of $(\ref{condXM})$ is equal to $I+J$ with 
\begin{eqnarray*}
I &=&E(f(X_{T_{k}}^{M})1_{\{U_{k}\geq\gamma (Z_{k},X^M_{T_{k}-})\}}\mid
X^M_{T_{k}-}=x)\quad and \\
J &=&E(f(X_{T_{k}}^{M})1_{\{U_{k}< \gamma (Z_{k},X_{T_{k}-}^{M})\}}\mid
X_{T_{k}-}^{M}=x).
\end{eqnarray*}
A simple calculation leads to
\begin{equation*}
I=f(x)P(U_{k} \geq \gamma (Z_{k},x))=f(x)\theta _{M,\gamma }(x)=\int_{\left\vert
z\right\vert >M+1}f(x+c_{M}(z,x))q_{M}(z,x)dz
\end{equation*}
where the last equality results from the fact that $c_{M}(z,x)=0$ for $
\left\vert z\right\vert >M+1.$ Moreover one can easily see that $
J=\int_{\left\vert z\right\vert \leq M+1}f(x+c_{M}(z,x))q_{M}(z,x)dz$ and $(
\ref{condXM})$ is proved. 

From the relation (\ref{condXM}) we construct a process $(\overline{X}^M_t)$ equal in law to $(X_t^M)$ on the following way. 

We denote by $\Psi _{t}(x)$
the solution of $\Psi _{t}(x)=x+\int_{0}^{t}g(\Psi _{s}(x))ds.$ We assume
that the times $T_{k},k\in \mathbb{N}$ are fixed and we consider a sequence  $
(z_{k})_{k\in \mathbb{N}}$ with $z_{k}\in \mathbb{R}^{d}.$ Then we define $x_{t},t\geq 0$ by $
x_{0}=x$ and, if $x_{T_{k}}$ is given, then 
\begin{eqnarray*}
x_{t} &=&\Psi _{t-T_{k}}(x_{T_{k}})\quad T_{k}\leq t<T_{k+1}, \\
x_{T_{k+1}} &=&x_{T_{k+1}^-}+c_{M}(z_{k+1},x_{T_{k+1}^-}).
\end{eqnarray*}
We remark that for $T_{k}\leq t<T_{k+1},x_{t}$ is a function of $
z_{1},...,z_{k}.$ Notice also that $x_{t}$ solves the equation
\begin{equation*}
x_{t}=x+\sum_{k=1}^{J_{t}^{M}}c_{M}(z_{k},x_{T_{k}^-})+\int_{0}^{t}g(x_{s})ds.
\end{equation*}
We consider now a sequence of random variables $(\overline{Z}_{k}),k\in
\mathbb{N}^{\ast }$ and we denote $\mathcal{G}_{k}=\sigma (T_{p},p\in \mathbb{N})\vee \sigma (
\overline{Z}_{p},p\leq k)$ and $\overline{X}_{t}^{M}=x_{t}(\overline{Z}
_{1},...,\overline{Z}_{J_{t}^{M}}).$ We assume that the law of $
\overline{Z}_{k+1}$ conditionally on $\mathcal{G}_{k}$ is given by 
\begin{equation*}
P(\overline{Z}_{k+1}\in dz\mid \mathcal{G}_{k})=q_{M}(x_{T_{k+1}^-}(\overline{
Z}_{1},...,\overline{Z}_{k}),z)dz=q_{M}(\overline{X}_{T_{k+1}^-}^{M},z)dz.
\end{equation*}
Clearly $\overline{X}_{t}^{M}$ satisfies the equation 
\begin{equation}
\overline{X}_{t}^{M}=x+\sum_{k=1}^{J_{t}^{M}}c_{M}(\overline{Z}_{k},\overline{X}
_{T_{k}-}^{M})+\int_{0}^{t}g(\overline{X}_{s}^{M})ds  \label{eq4}
\end{equation}
and $\overline{X}_{t}^{M}$ has the same law as $X_{t}^{M}.$
Moreover we can prove a little bit more. 

\begin{lem}\label{egaloi}
For a locally bounded and measurable
function $\psi :\mathbb{R}^{d}\rightarrow \mathbb{R}$ let 
\begin{equation*}
\overline{S}_{t}(\psi )=\sum_{k=1}^{J_{t}^M}(\Phi _{M}\psi )(\overline{Z}
_{k}),\quad S_{t}(\psi )=\sum_{k=1}^{J_{t}^M}(\Phi _{M}\psi
)(Z_{k})1_{\{\gamma (Z_{k},X^{M}(T_{k}-))> U_{k}\}},
\end{equation*}
then $(\overline{X}_{t}^{M},\overline{S}_{t}(\psi ))_{t\geq 0}$ has the same law
as $(X_{t}^{M},S_{t}(\psi ))_{t\geq 0}.$
\end{lem}

\begin{pf}
 Observing that  $(\overline{X}_{t}^{M},\overline{S}_{t}(\psi ))_{t\geq
0} $ solves a system of equations similar to $(\ref{eq4})$ but in dimension $
d+1$, it suffices to prove that $(\overline{X}_{t}^{M})_{t\geq 0}$ has the
same law as $(X_{t}^{M})_{t\geq 0}.$ This readily follows from 
\begin{equation*}
E(f(X_{T_{k+1}}^{M})\mid X_{T_{k+1}-}^{M}=x)=E(f(\overline{X}
_{T_{k+1}}^{M})\mid \overline{X}_{T_{k+1}-}^{M}=x)
\end{equation*}
which is a consequence of (\ref{condXM}). 

\end{pf}

\begin{rem}
Looking at the infinitesimal generator $L$ of $X$ it is clear that the
natural approximation of $X_{t}$ is $\overline{X}_{t}^{M}$ instead of $
X_{t}^{M}.$ But we use the representation given by $X_{t}^{M}$ for two
reasons. First  it is easier to obtain estimates for this process
because we have a stochastic equation and so we may use the stochastic
calculus associated to a Poisson point measure. Moreover, having this
equation in mind, gives a clear idea about the link with other approaches by
Malliavin calculus to the solution of a stochastic equation with jumps: we
mainly think to [B.G.J]. Remark that $X_t$ is solution of an equation with discontinuous coefficients so the approach developped by [B.G.J] does not
work. And if we consider the equation of $\overline{X}_{t}^{M}$ then the
underlying point measure depends on the solution of the equation so it is no
more a Poisson point measure.
\end{rem}

\subsection{The integration by parts formula}

The random variable $\overline{X}^M_t$ constructed previously is a simple functional but unfortunately its Malliavin
covariance matrix is degenerated. To avoid this problem we use a classical regularization procedure. Instead of the variable $\overline{X}^M_t$,
we consider the regularized one $F_M$ defined by
\begin{equation}
F_{M}=\overline{X}_{t}^{M}+\sqrt{U_{M}(t)}\times \Delta , \label{FM}
\end{equation}
where $\Delta$ is a $d-$dimensional standard gaussian variable independent of the variables $(\overline{Z}_k)_{k \geq 1}$ and $(T_k)_{k \geq 1}$ and $U_M(t)$ is defined by
\begin{eqnarray}
U_M(t) = t \int_{B_{M-1}^c} \underline{c}^2(z) \underline{\gamma}(z) d\mu(z). \label{UM}
\end{eqnarray}
We observe that  $F_M \in \mathcal{S}^d$ where $\mathcal{S}$ is  the space of simple functionals for the differential calculus  based on the variables $(\overline{Z}_k)_{k \in \mathbb{N}}$ with $\overline{Z}_0=(\Delta^r)_{1 \leq r \leq d}$ and $\overline{Z}_k=(\overline{Z}_k^r)_{1 \leq r \leq d}$ and
we are now in the framework of section $2$ by taking $\mathcal{G}=\sigma
(T_{k},k\in \mathbb{N})$ and defining the weights $(\pi_k)$ by $\pi_0^r=1$  and
$\pi_k^r=\Phi_M(\overline{Z}_k)$ for $1 \leq r \leq d$.
Conditionally on $\mathcal{G}$, the density of the law of $(\overline{Z}_{1},...,\overline{Z}_{J_{t}^M})$ is
given by 
\begin{equation*}
p_{M}(\omega ,z_{1},...,z_{J_{t}^M})=\prod_{j=1}^{J_{t}^M}q_{M}(z_{j},\Psi
_{T_{j}-T_{j-1}}(\overline{X}_{T_{j-1}}^{M}))
\end{equation*}
where $\overline{X}_{T_{j-1}}^{M}$ is a function of $z_{i}, 1 \leq i \leq j-1.$ We can check that $p_M$ satisfies the hypothesis H1 of section 2.

To clarify the notation, the derivative operator can be written in this framework for $F \in \mathcal{S}$ by $DF=(D_{k,r} F)$ where $D_{k,r}=\pi_k^r \partial_{\overline{Z}_k^r}$ for $k \geq 0$ and
$1 \leq r \leq d$. Consequently we deduce that  $D_{k,r} F_M^{r'}=D_{k,r} \overline{X}^{M,r'}_t$, for $k \geq 1$ and $D_{0,r} F_M^{r'}=\sqrt{U_M(t)} \delta_{r,r'}$ with $\delta_{r,r'}=0$ if $r \neq r'$,  $\delta_{r,r'}=1$  otherwise.

The Malliavin covariance matrix of $\overline{X}_{t}^{M}$ is equal to
\begin{equation*}
\sigma(\overline{X}_t^M)^{i,j}=\sum_{k=1}^{J_{t}^M}\sum_{r=1}^{d}D_{k,r}\overline{X}_{t}^{M,i}D_{k,r}\overline{X}_{t}^{M,j}
\end{equation*}
for $1 \leq i,j \leq d$ and finally the Malliavin covariance matrix of $F_M$ is given by
$$
\sigma(F_M)=\sigma(\overline{X}_t^M) + U_M(t) \times Id.
$$
Using the results of section $2$, we can state an integration by part formula and give a bound for the weight $H^q(F_M,1)$ in terms of the Sobolev norms of $F_M$, the divergence $L F_M$ and the determinant of  the inverse of the Malliavin covariance matrix $\det \sigma(F_M)$. The control of these last three quantities is rather technical and is studied in detail in section 4.

\begin{prop}\label{pIPPFM}
Assume hypotheses $3.0.$ $3.1.$ $3.2.$ and let $\phi :\mathbb{R}^{d}\rightarrow \mathbb{R}$ be a bounded smooth function with bounded derivatives.  For every multi-index $\beta=(\beta_1, \ldots \beta_q) \in \{1, \ldots,d\}^q$ such that $4d(3q-1)/t < \theta$

a) if  $3.3.a$ holds then
\begin{equation}
\left\vert E(\partial _{\beta }\phi (F_{M}))\right\vert \leq C_q \left\Vert \phi \right\Vert _{\infty }.  \label{IPPFMa}
\end{equation}

b) if $3.3.b$ holds then
\begin{equation}
\left\vert E(\partial _{\beta }\phi (F_{M}))\right\vert \leq C_q \left\Vert \phi \right\Vert _{\infty }
(1+\mu(B_{M+1})^q),  \label{IPPFMb}
\end{equation}
\end{prop}
{\bf Remark :} if $\theta= \infty$ then $\forall t>0$, we have an integration by parts formula for any order of derivation $q$. Conversely if $\theta$ is finite, we need to have $t$ large enough to integrate $q$ times by part.

\begin{pf}
The integration by parts formula $(\ref{IPP2})$ gives, for every smooth $\phi
:\mathbb{R}^{d}\rightarrow \mathbb{ R}$ and every multi-index $\beta =(\beta _{1},...,\beta
_{q})$
\begin{equation*}
E(\partial _{\beta }\phi (F_{M}))=E(\phi (F_{M})H_{\beta}^q(F_{M},1)),
\end{equation*}
and consequently
\begin{eqnarray*}
\left\vert E(\partial _{\beta}\phi (F_{M}))\right\vert  &\leq &\left\Vert
\phi \right\Vert _{\infty }E(\vert H_{\beta}^q(F_{M},1)\vert ) .
\end{eqnarray*}
So we just have to bound $\vert H_{\beta}^q(F_{M},1)\vert $. From the second part of Theorem \ref{thHq}  we have
$$
\vert H^q(F_M,1) \vert \leq C_q \frac{1}{ \vert \det \sigma(F_M) \vert^{3q-1}}(1+ \vert F_M\vert_{q+1}^{(6d+1)q} )( 1 + \vert L F_M \vert_{q-1}^q).
$$
Now from Lemma \ref{LFM} (see section 4), we have :

a) assuming $3.3.a$, for $l,p \geq 1$, 
$$
E \vert L F_M \vert_{l}^p  \leq C_{l,p};
$$

b) assuming $3.3.b$, for $l,p \geq 1$, 
$$
E \vert L F_M \vert_{l}^p \leq C_{l,p}(1+ \mu(B_{M+1})^p).
$$
Hence from  Lemma \ref{boundFM}, for $l,p \geq 1$
$$
E \vert F_M \vert_l^p \leq C_{l,p};
$$
and from Lemma \ref{lemsigmaFM} , we have for $p \geq 1$, $t>0$ such that $2dp/t<\theta$
$$
E \frac{1}{\det \sigma(F_M))^p} \leq C_p.
$$
The final result is then a straightforward consequence of Cauchy-Schwarz inequality. 
\end{pf}

\subsection{Estimates for the Fourier transform of $X_t$}
In this section, we prove Proposition \ref{fourierX}.

\begin{pf}
The proof consists first to approximate $X_t$ by $\overline{X}_t^M$ and then to apply the integration by parts formula.

{\bf Approximation.}
We have
\begin{equation*}
\left\vert E(e^{i \left\langle \xi ,X_{t} \right\rangle})\right\vert \leq \left\vert \xi \right\vert
E\left\vert X_{t}-\overline{X}_{t}^M\right\vert +\left\vert
E(e^{i \left\langle \xi ,\overline{X}_t^M \right\rangle}-e^{i \left\langle \xi ,F_{M} \right\rangle})\right\vert
+\left\vert E(e^{i \left\langle \xi, F_{M}\right\rangle})\right\vert .
\end{equation*}
From  $(\ref{approx})$ we deduce
\begin{equation*}
E(\left\vert X_{t}-\overline{X}_t^M \right\vert )\leq
\varepsilon _{M}=  t e^{Ct} \int_{B_{M}^{c}}\overline{c}(z)\overline{
\gamma }(z)d\mu (z).
\end{equation*}
Moreover
\begin{eqnarray*}
E(e^{i \left\langle \xi, \overline{X}_t^M \right\rangle}-e^{i \left\langle \xi ,F_{M}\right\rangle}) &=&E(e^{i \left\langle \xi ,
\overline{X}_t^M \right\rangle}(1-e^{i \left\langle \xi ,\sqrt{U_{M}(t)}\Delta \right\rangle}))
=E(e^{i \left\langle \xi ,\overline{X}_t^M \right\rangle})(1-e^{-\frac{1}{2}
\left\vert \xi \right\vert ^{2}U_{M}(t )}),
\end{eqnarray*}
so that 
\begin{equation*}
\left\vert E(e^{i \left\langle \xi, \overline{X}_t^M \right\rangle}-e^{i \left\langle \xi,
F_{M} \right\rangle})\right\vert \leq U_M(t) \frac{1}{2}\left\vert \xi
\right\vert ^{2}.
\end{equation*}
We conclude that
\begin{equation*}
\left\vert E(e^{i \left\langle \xi ,X_{t} \right\rangle})\right\vert \leq U_M(t)
\frac{1}{2}\left\vert \xi \right\vert ^{2}+\left\vert \xi \right\vert te^{Ct}
\int_{B_{M}^{c}}\overline{c}(z)d\mu (z))+\left\vert E(e^{i \left\langle \xi,
F_{M} \right\rangle})\right\vert .
\end{equation*}
\textbf{Integration by parts.} We denote $e_{\xi }(x)=\exp (i \left\langle \xi ,x \right\rangle)$
and we have $\partial _{\beta}e_{\xi }(x)=i^{\left\vert \beta \right\vert }\xi_{\beta_1} \ldots \xi_{\beta_q}
e_{\xi }(x)$. Consequently 

a) assuming $3.3.a$ and applying $(\ref{IPPFMa})$ for $\beta $ such that $\left\vert \beta
\right\vert =q$ we obtain
\begin{equation*}
\left\vert E(e^{i\ \left\langle \xi ,F_{M} \right\rangle})\right\vert  \leq 
\frac{C_q}{\vert \xi \vert^q},
\end{equation*}

b) assuming $3.3.b$, we obtain similarly from $(\ref{IPPFMb})$
\begin{equation*}
\vert \xi_{\beta_1} \ldots \xi_{\beta_q} \vert \left\vert E(e^{i\ \left\langle \xi ,F_{M} \right\rangle})\right\vert  =
\left\vert E(\partial _{\beta}e_{\xi }(F_{M}))\right\vert \leq C_q(1+\mu(B_{M+1})^q),
\end{equation*}
and then
\begin{equation*}
\left\vert E(e^{i\ \left\langle \xi ,F_{M} \right\rangle})\right\vert  \leq 
\frac{C_q}{\vert \xi \vert^q}(1+\mu(B_{M+1})^q),
\end{equation*}

and the proposition is proved.

\end{pf}

\section{Sobolev norms-Divergence-Covariance matrix}
\subsection{Sobolev norms}
We prove in this section that $\forall l \geq 1$ and $\forall p \geq 1$ $E\vert F_M\vert_l^p \leq C_{l,p}$. We begin this section with a preliminary lemma which will be also useful to control the covariance matrix.

\subsubsection{Preliminary}

We consider a Poisson point measure 
$N(ds,dz,du)$ on $\mathbb{R}^{d}\times \mathbb{R}_{+}$ with compensator $\mu (dz)\times
1_{(0,\infty )}(u)du$ and  two non negative measurable
functions $f,g:\mathbb{R}^{d}\rightarrow \mathbb{R}_{+}$. For a measurable set $B\subset \mathbb{R}^{d}$
we denote $B_{g}=\{(z,u):z\in B,u<g(z)\}\subset \mathbb{R}^{d}\times \mathbb{R}_{+} $ and we
consider the process
\begin{equation*}
N_{t}(1_{B_{g}}f):=\int_{0}^{t}\int_{B_{g}}f(z)N(ds,dz,du).
\end{equation*}
Moreover we note $\nu _{g}(dz)=g(z)d\mu (z)$ and 
\begin{equation*}
\alpha _{g,f}(s)=\int_{\mathbb{R}^{d}}(1-e^{-sf(z)})d\nu _{g}(dz),\quad \beta
_{B,g,f}(s)=\int_{B^{c}}(1-e^{-sf(z)})d\nu _{g}(dz).
\end{equation*}
We have the following result.
\begin{lem}\label{laplace}
Let $\phi(s)=Ee^{-s N_t(f 1_{B_g})}$ the Laplace transform of the random variable $ N_t(f 1_{B_g})$ then we have
\begin{eqnarray*}
\phi(s)=e^{-t(\alpha _{g,f}(s)-\beta_{B,g,f}(s))}.
\end{eqnarray*}
\end{lem}

\begin{pf}
From It\^{o}'s formula we have 
\begin{equation*}
\exp (-sN_{t}(f1_{B_{g}}))=1-\int_{0}^{t}\int_{\mathbb{R}^{d}\times \mathbb{R}_{+}}\exp
(-s(N_{r-}(f1_{B_{g}})))(1-\exp (-sf(z)1_{B_{g}}(z,u)))dN(r,z,u)
\end{equation*}
and consequently 
\begin{equation*}
E(\exp (-sN_{t}(f1_{B_{g}})))=1-\int_{0}^{t}E(\exp
(-s(N_{r-}(f1_{B_{g}}))\int_{\mathbb{R}^{d}\times \mathbb{R}_{+}}(1-\exp
(-sf(z)1_{B_{g}}(z,u)))d\mu (z)dudr.
\end{equation*}
But 
\begin{eqnarray*}
\int_{\mathbb{R}^{d}\times \mathbb{R}_{+}}(1-\exp (-sf(z)1_{B_{g}}(z,u)))d\mu (z)du
&=&\int_{\mathbb{R}^{d}\times \mathbb{R}_{+}}1_{B_{g}}(z,u)(1-\exp (-sf(z)))d\mu (z)du \\
&=&\int_{\mathbb{R}^{d}}1_{B}(z)(1-\exp (-sf(z)))\int_{\mathbb{R}_{+}}1_{\{u<
g(z)}\}dud\mu (z) \\
&=&\int_{B}(1-\exp (-sf(z)))g(z)d\mu (z)=\alpha _{g,f}(s)-\beta
_{B,g,f}(s),
\end{eqnarray*}
It follows that 
\begin{equation*}
E(\exp (-sN_{t}(f1_{B_{g}})))=\exp (-t(\alpha _{g,f}(s)-\beta
_{B,g,f}(s))).
\end{equation*}
\end{pf}

\subsubsection{Bound for $\vert \overline{X}_t^M\vert_l$}
In this section, we use the notation $\overline{c}_1(z) = \sup_x \vert \nabla_x c(z,x) \vert$. Under hypothesis $3.3.i$
we have $\overline{c}_1(z) \leq \overline{c}(z)$, but we introduce this notation to highlight the dependence on the first
derivative of the function $c$.
\begin{lem}\label{boundXM}
Let $(\overline{X}_t^M)$ the process solution of equation $(\ref{eq4})$ then under hypotheses $3.0.$, $3.1.$ and $3.2.$ we have
$\forall l \geq 1$,  
$$
\sup_{s \leq t}\vert \overline{X}_s^M\vert_{1,l} \leq C_l (1+ \sum_{k=1}^{J_t^M} \overline{c}(\overline{Z}_k ))^{l \times l!} \sup_{s \leq t} (\mathcal{E}_s^M)^{l\times l!}
$$
where $C_l$ is an universal constant and where $\mathcal{E}_t^M$  is solution of the linear equation
\begin{eqnarray}
\mathcal{E}_t^M=1 + C_l\sum_{k=1}^{J_t^M} \overline{c}_1(\overline{Z}_k ) \mathcal{E}_{T_k-}^M + C_l \int_0^t \mathcal{E}_s^M ds  .\label{exp}
\end{eqnarray}
Consequently
 $ \forall l,p \geq 1$
$$
\sup_M E\sup_{s \leq t}\vert \overline{X}_s^M\vert_{1,l}^p < \infty
$$
\end{lem}
Before proving this lemma we first give a result which is a straightforward consequence of lemma \ref{CR} and formula $(\ref{CRk})$.

\begin{lem}\label{CRnorml}
Let $\phi : \mathbb{R}^d \mapsto \mathbb{R}$ a $\mathcal{C}^{\infty}$ function and $F \in \mathcal{S}^d$ then $\forall l \geq 1$ we have
$$
\vert \phi(F) \vert_{1,l} \leq \vert \nabla \phi(F) \vert \vert F \vert_{1,l} + C_l \sup_{2 \leq \vert \beta \vert \leq l} \vert \partial_{\beta} \phi(F) \vert 
\vert F \vert_{1,l-1}^l.
$$
\end{lem}
We proceed now to the proof of Lemma \ref{boundXM}.
\begin{pf}
We first recall that from hypothesis $3.0.$, $g$ and its derivatives are bounded and from hypothesis $3.2.i)$ the coefficient $c$ as well as its derivatives are bounded by the function $\overline{c} $. Now the truncated coefficient $c_M$ of equation (\ref{eq4}) is equal to $c_M=c \times \phi_M$ where $\phi_M$ is a $\mathcal{C}^{\infty}$ bounded function with derivatives uniformly bounded with respect to $M$. Consequently using Lemma \ref{CRnorml} we obtain for $l \geq 1$
\begin{eqnarray*}
\vert \overline{X}_t^M \vert_{1,l} \leq C_l \left( A_{t,l-1} +  \sum_{k=1}^{J_t^M} \overline{c}_1(\overline{Z}_k )\vert \overline{X}_{T_k-}^M \vert_{1,l}
+ \int_0^t \vert \overline{X}_s^M \vert_{1,l} ds \right),
\end{eqnarray*}
with
$$
A_{t,l-1}=\sum_{k=1}^{J_t^M} \overline{c}(\overline{Z}_k )(\vert\overline{Z}_k \vert_{1,l} + \vert\overline{Z}_k \vert_{1,l-1}^l + 
\vert \overline{X}_{T_k-}^M \vert_{1,l-1}^l)
+ \int_0^t \vert \overline{X}_s^M \vert_{1,l-1}^l ds .
$$
This gives
\begin{eqnarray}
\forall s \leq t \quad \vert \overline{X}_s^M \vert_{1,l} \leq A_{t ,l-1} \mathcal{E}_s^M ,\label{BXM} 
\end{eqnarray}
Under hypotheses $3.0.$ $3.1.$ and $3.2.$ we have 
$$
\forall p \geq 1 \quad  E(\sup_{s \leq t}\vert \mathcal{E}_t^M \vert^p) \leq C_p.
$$
Now one can easily check that for $l \geq 1$
$$
\vert\overline{Z}_k \vert_{1,l} \leq \vert \pi_k \vert_{l-1},
$$
but since $\pi_k= \phi_M(\overline{Z}_k)$ we deduce from Lemma \ref{CRnorml} that
$$
\vert\overline{Z}_k \vert_{1,l} \leq 1 + C_l(\vert\overline{Z}_k \vert_{1,l-1} + \vert\overline{Z}_k \vert_{1,l-2}^{l-1}).
$$
Observing that $\vert\overline{Z}_k \vert_{1,1}= \vert D \overline{Z}_k \vert= \vert \pi_k \vert \leq 1$ we conclude that $\forall l \geq 1$
$$
\vert\overline{Z}_k \vert_{1,l} \leq C_l.
$$
This gives
\begin{eqnarray}
A_{t,l-1}\leq t C_l (1+ \sup_{s \leq t}\vert \overline{X}_s^M \vert_{1,l-1})^l(1+\sum_{k=1}^{J_t^M} \overline{c}(\overline{Z}_k )). \label{At}
\end{eqnarray}
From this inequality we can prove easily Lemma \ref{boundXM} by induction. For $l=1$ we remark that
\begin{eqnarray*}
\forall s \leq t \quad \vert \overline{X}_s^M \vert_{1,1} \leq A_{t,0} \mathcal{E}_s^M , \quad \mbox{with} \quad 
A_{t,0}=\sum_{k=1}^{J_t^M} \overline{c}(\overline{Z}_k ),
\end{eqnarray*}
and the result is true. To complete the proof of lemma \ref{boundXM}, we prove that $\forall p \geq 1$
$$
E\left(\sum_{k=1}^{J_t^M} \overline{c}(\overline{Z}_k )\right)^p \leq C_p.
$$
We have the equality in law
$$
\sum_{k=1}^{J_t^M} \overline{c}(\overline{Z}_k )\backsimeq \int_0^t  \int_E \overline{c}(z)1_{u<\gamma(z, X_{T_k-}^M)} 1_{B_{M+1}}(z) 1_{[0,2 \overline{C}]}(u) N(ds,dz,du),
$$  
moreover using the notations of section $4.1.1.$ we have
$$
\int_0^t  \int_E \overline{c}(z)1_{u<\gamma(z, X_{T_k-}^M)} 1_{B_{M+1}}(z) 1_{[0,2 \overline{C}]}(u) N(ds,dz,du) \leq N_t(1_{B_{\overline{\gamma}}} \overline{c})
$$
with $B_{\overline{\gamma}}= \{(z,u); z \in B_{M+1}; 0<u< \overline{\gamma}(z) \}$.
From Lemma \ref{laplace} it follows that
$$
Ee^{-sN_t(1_{B_{\overline{\gamma}}} \overline{c})} = \exp(-t \int_{B_{M+1}}(1-e^{-s \overline{c}(z)}) \overline{\gamma} (z)d \mu(z))
$$
and since from hypotheses $3.1.$ and $3.2.$,  $ \int_{\mathbb{R}^d} \vert \overline{c}(z) \overline{\gamma}(z) \vert d \mu(z) < \infty$ we deduce that
$\forall p \geq 1$, $EN_t(1_{B_{\overline{\gamma}}} \overline{c})^p = t^p (\int_{B_{M+1}} \vert \overline{c}(z) \overline{\gamma}(z) \vert d \mu(z))^p \leq C_p$ where the constant $C_p$ does not depend on $M$. This achieves the proof of 
Lemma \ref{boundXM}.
\end{pf}
\subsubsection{Bound for $\vert F_M \vert_l$}

\begin{lem}\label{boundFM}
Under hypotheses $3.0.$, $3.1.$ and $3.2.$ we have
$$
\forall l,p \geq 1 \quad E \vert F_M \vert_l^p \leq C_{l,p}.
$$
\end{lem}
We have $F_M =\overline{X}_t^M+\sqrt{U_M(t)} \Delta$ and then $\vert F_M \vert_l \leq \vert \overline{X}_t^M \vert_l+\sqrt{U_M(t)} \vert \Delta \vert_l$. But $\vert \Delta \vert_l \leq \vert \Delta \vert +d$ and $U_M(t) \leq t \int_{ \mathbb{R}^d} \overline{c}^2(z) \overline{ \gamma}(z) d \mu(z) < \infty$. So the conclusion of Lemma  \ref{boundFM} follows from Lemma \ref{boundXM}.

\subsection{Divergence}
In this section our goal is to bound $\vert L F_M \vert_l$ for $l \geq 0$. From the definition of the divergence operator $L$
we have $L F_M^r=L \overline{X}_t^{M,r}- \Delta^r$ and then
$$
\vert L F_M \vert_l \leq \vert  L \overline{X}_t^{M}\vert_l + \vert \Delta \vert +d,
$$
so we just have to bound $\vert  L \overline{X}_t^{M}\vert_l  $. We proceed as in the previous section and we first state a lemma
similar to Lemma \ref{CRnorml}.

\begin{lem}\label{Lnorml}
Let $\phi : \mathbb{R}^d \mapsto \mathbb{R}$ a $\mathcal{C}^{\infty}$ function and $F \in \mathcal{S}^d$ then $\forall l \geq 1$ we have
$$
\begin{array}{lll}
\vert L \phi(F) \vert_{l} & \leq & \vert \nabla \phi(F) \vert \vert L F \vert_{l} + C_l \sup_{2 \leq \vert \beta \vert \leq l+2} \vert \partial_{\beta} \phi(F) \vert 
(1+\vert F \vert_{l}^l)(\vert L F \vert_{l-1} + \vert F \vert_{1,l+1}^2), \\
& \leq & \vert \nabla \phi(F) \vert \vert L F \vert_{l} + C_l \sup_{2 \leq \vert \beta \vert \leq l+2} \vert \partial_{\beta} \phi(F) \vert 
(1+\vert F \vert_{l+1}^{l+2})(1 + \vert L F \vert_{l-1}).
\end{array}
$$
For $l=0$, we have
$$
\vert L \phi(F) \vert \leq  \nabla \phi(F) \vert \vert L F \vert +\sup_{ \beta =2} \vert \partial_{\beta} \phi{F} \vert \vert F \vert_{1,1}^2.
$$
\end{lem}
The proof follows from (\ref{CRL}) and Lemma \ref{CRnorml} and we omit it. 

Next we give a bound for $\vert L \overline{Z}_k \vert_l$.
We recall the notation 
\begin{eqnarray*}
 \overline{\gamma}_{\ln}^{z,l}(z)= \sup_x \sup_{1 \leq \vert \beta \vert \leq l} \vert \partial_{\beta,z} \ln \gamma(z,x) \vert, \quad
\overline{h}_{\ln}^{l}(z)=  \sup_{1 \leq \vert \beta \vert \leq l} \vert \partial_{\beta} \ln h(z) \vert, 
\quad \overline{\theta}_{\ln}^{l}= \sup_x  \sup_{1 \leq \vert \beta \vert \leq l} \vert \partial_{\beta} \ln \theta_{M, \gamma}(x) \vert, 
\end{eqnarray*}
\begin{eqnarray*} 
 \overline{\gamma}_{\ln}^{x,l}(z)= \sup_x \sup_{1 \leq \vert \beta \vert \leq l} \vert \partial_{\beta,x} \ln \gamma(z,x)\vert, 
 \quad  \overline{\gamma}^{x,l}= \sup_z \sup_x \sup_{1 \leq \vert \beta \vert \leq l} \vert \partial_{\beta,x}  \gamma(z,x)\vert.
\end{eqnarray*}
\begin{lem}\label{LZ}
Assuming hypotheses $3.0.$, $3.1.$, $3.2$  and $3.3.$, we have $\forall l \geq 0$ and  $ \forall k \leq J_t^M$
$$
\vert L \overline{Z}_k \vert_l \leq C_l ( \overline{\gamma}_{\ln}^{z,l+1}(\overline{Z}_k)+\overline{h}_{\ln}^{z,l+1}(\overline{Z}_k) + \sup_{s \leq t} \vert \overline{X}_s^M \vert_{l+1}^{l+1}  \sum_{j=k+1}^{J_t^M}  \overline{\theta}_{\ln}^{l+1}1_{B(z_M^*,1)}(\overline{Z}_j )+
\overline{\gamma}_{\ln}^{x,l+1}(\overline{Z}_j))),
$$
with $\overline{\theta}_{\ln}^{l} \leq C_l ( \overline{\gamma}^{x,l})^l.$

In addition, if we assume $3.3.a.$,  we obtain $\forall p \geq 1$
$$
E\sup_{k\leq J_t^M} \vert L \overline{Z}_k \vert_l ^p \leq C_{p,l}.
$$
On the other hand, assuming $3.3.b$, we have 
 $\forall p \geq 1$
$$
E\sup_{k\leq J_t^M} \vert L \overline{Z}_k \vert_l ^p \leq C_{p,l} (1+\mu(B_{M+1})^p)
$$
\end{lem}
\begin{pf}
We first recall that we have proved in the preceding section that $\forall l \geq 1$, $\vert \overline{Z}_k \vert_l \leq C_l$.
Now $L \overline{Z}_k^r= \delta(D \overline{Z}_k^r)$ and since $D_{k,r} \overline{Z}_k^r= \pi_k$ we obtain
$$
L \overline{Z}_k^r= -\partial_{k,r}( \pi_k^2)- \pi_k D_{k,r} \ln p_M,
$$
this leads to
$$
\vert L \overline{Z}_k^r \vert_l \leq C_l (1+ \vert D_{k,r} \ln p_M \vert_l).
$$
Recalling that $\ln p_M=\sum_{j=1}^{J_t^M} \ln q_M(\overline{Z}_j, \overline{X}^M_{T_j-})$ and that $ \overline{X}^M_{T_j-}$
depends on $\overline{Z}_k$ for $k \leq j-1$ we obtain
$$
D_{k,r} \ln p_M = D_{k,r} \ln q_M(\overline{Z}_k, \overline{X}^M_{T_k-}) +\sum_{j=k+1}^{J_t^M}D_{k,r} \ln q_M(\overline{Z}_j, \overline{X}^M_{T_j-})
$$
But on $\{ \pi_k>0 \}$, we have $q_M(\overline{Z}_k, \overline{X}^M_{T_k-}) =C \gamma(\overline{Z}_k, \overline{X}^M_{T_k-}) 
h(\overline{Z}_k)$, and then
$$
D_{k,r} \ln q_M(\overline{Z}_k, \overline{X}^M_{T_k-}) = D_{k,r} \ln \gamma (\overline{Z}_k, \overline{X}^M_{T_k-}) +D_{k,r} \ln h(\overline{Z}_k).
$$
Now for $j \geq k+1$, if $\vert\overline{Z}_j-z_M^* \vert <1$ then
$$
\ln q_M(\overline{Z}_j, \overline{X}^M_{T_j-})=\ln \varphi(\overline{Z}_j-z_M^*  ) + \ln \theta_{M, \gamma}( \overline{X}^M_{T_j-})
$$
consequently
$$
D_{k,r}\ln q_M(\overline{Z}_j, \overline{X}^M_{T_j-})=D_{k,r}\ln \theta_{M, \gamma}( \overline{X}^M_{T_j-}),
$$
and if $\overline{Z}_j \in B_{M+1}$ then
$$
D_{k,r}\ln q_M(\overline{Z}_j, \overline{X}^M_{T_j-})=D_{k,r}\ln \gamma( \overline{Z}_j,\overline{X}^M_{T_j-})
$$
and finally
$$
D_{k,r}\ln q_M(\overline{Z}_j, \overline{X}^M_{T_j-})=D_{k,r}\ln \theta_{M, \gamma}( \overline{X}^M_{T_j-})1_{B(z_M^*,1)}( \overline{Z}_j)
+ D_{k,r}\ln \gamma( \overline{Z}_j,\overline{X}^M_{T_j-})1_{B_{M+1}}(\overline{Z}_j ).
$$
It is worth to note that this random variable is a simple variable as defined in section 2.

Putting this together, it yields
\begin{eqnarray*}
\vert D_{k,r} \ln p_M \vert_l  & \leq & \vert D_{k,r} \ln \gamma (\overline{Z}_k, \overline{X}^M_{T_k-}) \vert_l+ \vert D_{k,r} \ln h(\overline{Z}_k)
\vert_l \\ 
& &+ \sum_{j=k+1}^{J_t^M} (\vert D_{k,r}\ln \theta_{M, \gamma}( \overline{X}^M_{T_j-}) 1_{B(z_M^*,1)}( \overline{Z}_j)
\vert_l+
\vert D_{k,r}\ln \gamma( \overline{Z}_j,\overline{X}^M_{T_j-})\vert_l ).
\end{eqnarray*}
Applying Lemma \ref{CRnorml}, this gives
\begin{eqnarray*}
\vert D_{k,r} \ln p_M \vert_l \leq  ( \overline{\gamma}_{\ln}^{z,l+1}(\overline{Z}_k)+\overline{h}_{\ln}^{z,l+1}(\overline{Z}_k))
\vert \overline{Z}_k\vert_{1,l+1}^{l+1} 
+  \sum_{j=k+1}^{J_t^M}(  \overline{\theta}_{\ln}^{l+1}1_{B(z_M^*,1)}( \overline{Z}_j)
+
\overline{\gamma}_{\ln}^{x,l+1}(\overline{Z}_j)) \vert \overline{X}^M_{T_j-})\vert_{1,l+1}^{l+1} .
\end{eqnarray*} 
We obtain then, for $k \leq J_t^M$
$$
\vert L \overline{Z}_k \vert_l  \leq  C_l (\overline{\gamma}_{\ln}^{z,l+1}(\overline{Z}_k)+\overline{h}_{\ln}^{z,l+1}(\overline{Z}_k) + \sup_{s \leq t} \vert \overline{X}_s^M \vert_{l+1}^{l+1}  \sum_{j=k+1}^{J_t^M} ( \overline{\theta}_{\ln}^{l+1}1_{B(z_M^*,1)}( \overline{Z}_j) 
 + \overline{\gamma}_{\ln}^{x,l+1}(\overline{Z}_j))).
$$
Now from the definition of $\theta_{M, \gamma}$, we have
$$
\partial_{\beta}\theta_{M, \gamma}(x) =-\frac{1}{2 \overline{C} \mu(B_{M+1})} \int_{B_{M+1} }\partial_{\beta,x} \gamma(z,x) d \mu(z).
$$
Then assuming $3.3.$ and recalling that $1/2 \leq \theta_{M, \gamma}(x) \leq 1$, we obtain
$$
\overline{\theta}_{\ln}^l \leq C_l (\overline{\gamma}^{x,l})^l
$$
this finally gives
$$
\vert L \overline{Z}_k \vert_l  \leq  C_l (\overline{\gamma}_{\ln}^{z,l+1}(\overline{Z}_k)+\overline{h}_{\ln}^{z,l+1}(\overline{Z}_k) + \sup_{s \leq t} \vert \overline{X}_s^M \vert_{l+1}^{l+1}  \sum_{j=k+1}^{J_t^M} ( (\overline{\gamma}^{x,l+1})^{l+1}1_{B(z_M^*,1)}( \overline{Z}_j) 
 + \overline{\gamma}_{\ln}^{x,l+1}(\overline{Z}_j))).
$$
The first part of Lemma \ref{LZ} is proved. Moreover, we can check that from $3.3$, we have $\forall p \geq 1$
$$
E(\sum_{j=1}^{J_t^M} 1_{B(z_M^*,1)}( \overline{Z}_j))^p \leq t^p \sup_{z^*} (\int_{B(z^*,1)} \overline{\gamma}(z) d \mu(z))^p < \infty.
$$
Now assuming $3.3.a$, we have $\forall p \geq 1$
$$
E(\sum_{j=1}^{J_t^M} \overline{\gamma}_{\ln}^{x,l+1}(\overline{Z}_j))^p \leq t^p (\int \overline{\gamma}_{\ln}^{x,l+1}(z)  \overline{\gamma}(z) d \mu(z))^p < \infty,
$$
then the second part of Lemma \ref{LZ} follows from Lemma \ref{boundXM} and Cauchy-Schwarz inequality.
At last, assuming $3.3.b$, 
we check that $\sum_{j=1}^{J_t^M} \overline{\gamma}_{\ln}^{x,l+1}(\overline{Z}_j) \leq  \overline{\gamma}_{\ln}^{x,l+1} J_t^M$, and the third part follows easily. 
\end{pf}

We can now state the main lemma of this section.
\begin{lem}\label{LX}
Assuming hypotheses $3.0.$, $3.1.$  and $3.2.$, we have $\forall l \geq 0$ 
$$
\sup_{s \leq t} \vert L \overline{X}_s^M  \vert_l \leq B_{t,l}^M (1+ \sup_{k \leq J_t^M}  \vert L \overline{Z}_k \vert_l) ,
$$
where $B_{t,l}^M$ is a random variable such that $\forall p \geq 1 $, $E (B_{t,l}^M)^p \leq C_p$ for a constant $C_p$ independent
on $M$. More precisely we have
$$
B_{t,l}^M \leq C_l (1+ \sum_{k=1}^{J_t^M} \overline{c}(\overline{Z}_k))^{l+1} (1+ \sup_{s \leq t} \vert \overline{X}_s^M \vert_{l+1}^{l+2})^{l +1}\sup_{s \leq t} (\mathcal{E}_s^M)^{l+1},
$$
where $\mathcal{E}_s$ is solution of (\ref{exp}).
\end{lem}
\begin{pf}
We proceed by induction. From equation (\ref{eq4}) we have
$$
L \overline{X}_t^M=\sum_{k=1}^{J_t^M} L c_M(\overline{Z}_k, \overline{X}^M_{T_k-}) + \int_0^t L g(\overline{X}_s^M ) ds.
$$
For $l=0$, the second part of Lemma \ref{Lnorml} gives
$$
\vert L \overline{X}_t^M \vert \leq B_{t,0}+C \left(\sum_{k=1}^{J_t^M} \overline{c}_1(\overline{Z}_k)\vert L \overline{X}^M_{T_k-}\vert + \int_0^t \vert L \overline{X}_s^M \vert ds \right)
$$
with
$$
B_{t,0}=C \left( \sum_{k=1}^{J_t^M} \overline{c}(\overline{Z}_k)(\vert L \overline{Z}_k \vert +\vert \overline{Z}_k \vert _{1,1}^2+ \vert  \overline{X}^M_{T_k-}\vert_{1,1}^2) +  \int_0^t \vert  \overline{X}_s^M \vert _{1,1}^2 ds\right).
$$
This gives
$$
\forall s \leq t , \quad \vert L \overline{X}_s^M \vert \leq B_{t,0}\mathcal{E}_s^M,
$$
where $\mathcal{E}_s^M$ is solution of (\ref{exp}) and 
$$
B_{t,0} \leq C (1+ \sum_{k=1}^{J_t^M} \overline{c}(\overline{Z}_k)) (1+ \sup_{s \leq t} \vert \overline{X}_s^M \vert_{1}^{2})
(1+\sup_{k \leq J_t^M}  \vert L \overline{Z}_k \vert) .
$$
Consequently Lemma \ref{LX} is proved for $l=0$.

For $l >0$, we obtain similarly from Lemma \ref{Lnorml}
$$
\vert L \overline{X}_t^M \vert_l \leq B_{t,l-1}+C_l \left(\sum_{k=1}^{J_t^M} \overline{c}_1(\overline{Z}_k)\vert L \overline{X}^M_{T_k-}\vert_l + \int_0^t \vert L \overline{X}_s^M \vert_l ds \right)
$$
with
$$
\begin{array}{lll}
B_{t,l-1} &=&C_l  \sum_{k=1}^{J_t^M} \overline{c}(\overline{Z}_k)(\vert L \overline{Z}_k \vert_l +1+\vert L \overline{X}^M_{T_k-}\vert_{l-1})(1+
\vert \overline{Z}_k \vert _{l+1}^{l+2}+ \vert  \overline{X}^M_{T_k-}\vert_{l+1}^{l+2}) \\
& &+ C_l \int_0^t (1+\vert L\overline{X}^M_{T_k-} \vert_{l-1})(1+\vert  \overline{X}_s^M \vert _{l+1}^{l+2}) ds.
\end{array}
$$
We deduce then that
$$
\begin{array}{lll}
B_{t,l-1}& \leq &C_l  (1+ \sup_{s \leq t} \vert L \overline{X}_s^M \vert_{l -1})(1+\sup_{s \leq t} \vert  \overline{X}_s^M \vert_{l+1}^{l+2} ) (1+\sum_{k=1}^{J_t^M} \overline{c}(\overline{Z}_k) )\\
 & & +C_l \sup_{k \leq J_t^M}  \vert L \overline{Z}_k \vert_l(1+\sup_{s \leq t} \vert  \overline{X}_s^M \vert_{l+1}^{l+2}) \sum_{k=1}^{J_t^M} \overline{c}(\overline{Z}_k),
\end{array}
$$
now from the induction hypothesis, we have
$$
\begin{array}{lll}
B_{t,l-1}& \leq &C_l  (1+\sup_{s \leq t} \vert  \overline{X}_s^M \vert_{l+1}^{l+2} )^{l+1} (1+\sum_{k=1}^{J_t^M} \overline{c}(\overline{Z}_k) )^{l+1} \sup_{s \leq t} (\mathcal{E}_s^M)^{l}(1+ \sup_{k \leq J_t^M}  \vert L \overline{Z}_k \vert_{l-1})\\
 & & +C_l \sup_{k \leq J_t^M}  \vert L \overline{Z}_k \vert_l(1+\sup_{s \leq t} \vert  \overline{X}_s^M \vert_{l+1}^{l+2}) \sum_{k=1}^{J_t^M} \overline{c}(\overline{Z}_k),
\end{array}
$$
this leads to
$$
\forall s \leq t \quad \vert L \overline{X}_s^M \vert_l \leq B_{t,l}^M(1+ \sup_{k \leq J_t^M}  \vert L \overline{Z}_k \vert_{l}),
$$
with
$$
B_{t,l}^M \leq C_l  (1+\sup_{s \leq t} \vert  \overline{X}_s^M \vert_{l+1}^{l+2} )^{l+1} (1+\sum_{k=1}^{J_t^M} \overline{c}(\overline{Z}_k) )^{l+1} \sup_{s \leq t} (\mathcal{E}_s^M)^{l+1}.
$$
From Lemma \ref{boundXM}, we observe that $E  (B_{t,l}^M)^p < C_p$.

\end{pf}

Finally recalling that 
$$
\vert L F_M \vert_l \leq \vert L \overline{X}^M_t \vert_l + \vert \Delta \vert + d
$$
and combining Lemma \ref{boundXM}, Lemma \ref{LZ} and Lemma \ref{LX} we deduce easily the following lemma.
\begin{lem} \label{LFM}
Assuming hypotheses $3.0.$, $3.1.$ and $3.2.$, we have $\forall l,p \geq 1$

a) if $3.3.a$ holds, $E \vert L F_M \vert_l^p  \leq C_{l,p}$;

b) if  $3.3.b$ holds, $E \vert L F_M \vert_l^p \leq C_{l,p}(1+ \mu(B_{M+1})^p)$.
\end{lem}

\subsection{The covariance matrix}

\subsubsection{Preliminaries}

We consider an abstract measurable space $E,$ a measure $\nu $ on this space
and a non negative measurable function $f:E\rightarrow \mathbb{R}_{+}$ such that $
\int fd\nu <\infty .$ For $t>0$ and $p \geq 1$ we note
\begin{equation*}
\alpha _{f}(t)=\int_E (1-e^{-tf(a)})d\nu (a)\quad and\quad
I_{t}^{p}(f)=\int_{0}^{\infty }s^{p-1}e^{-t\alpha _{f}(s)}ds.
\end{equation*}

\begin{lem}\label{prelem}
i) Suppose that for  $p\geq 1$ and $t>0$
\begin{equation}
\underline{\lim }_{u\rightarrow \infty }\frac{1}{\ln u} \alpha_f(u) > p/t   \label{A1}
\end{equation}
then   $I_{t}^{p}(f)<\infty $. 

ii) A sufficient condition for $(\ref{A1})$ is
\begin{equation}
\underline{\lim }_{u\rightarrow \infty }\frac{1}{\ln u}\nu (f\geq \frac{1
}{u})> p/t  .\label{A2}
\end{equation}
In particular, if $\underline{\lim }_{u\rightarrow \infty }\frac{1}{\ln u}\nu (f\geq \frac{1
}{u})= \infty$ then $\forall p \geq 1$ and $\forall t>0$, $I_t^p(f)< + \infty$.
\end{lem}
We remark that if $\nu$ is finite then (\ref{A2}) can not be satisfied.
\begin{pf} i) 
From (\ref{A1}) one can find $\varepsilon>0$ such that as $s$ goes to infinity $s^{p-1} e^{-t \alpha_f(s)} \leq 1/ s^{1+ \varepsilon}$
and consequently $I_{t}^{p}(f)<\infty $.

ii) With the notation $n(dz)=\nu
\circ f^{-1}(dz)$ we have 
\begin{equation*}
\alpha _{f}(u)=\int_{0}^{\infty }(1-e^{-uz})dn(z)=\int_{0}^{\infty }e^{-y}n(
\frac{y}{u},\infty )dy.
\end{equation*}
Using Fatou's lemma and (\ref{A2}) we obtain
\begin{eqnarray*}
\underline{\lim }_{u\rightarrow \infty }\frac{1}{\ln u}\int_{0}^{\infty }e^{-y}n(
\frac{y}{u},\infty )dy &\geq &\int_{0}^{\infty }e^{-y}\underline{\lim }
_{u\rightarrow \infty }\frac{1}{ \ln u }n(\frac{y}{u},\infty )dy >p/t.
\end{eqnarray*}

\end{pf}

We come now back to the framework of section $4.1.1$ and we consider the Poisson point measure 
$N(ds,dz,du)$ on $\mathbb{R}^{d}\times \mathbb{R}_{+}$ with compensator $\mu (dz)\times
1_{(0,\infty )}(u)du.$ 
We recall that
\begin{equation*}
N_{t}(1_{B_{g}}f):=\int_{0}^{t}\int_{B_{g}}f(z)N(ds,dz,du),
\end{equation*}
for $f,g:\mathbb{R}^{d}\rightarrow \mathbb{R}_{+}$ and $B_{g}=\{(z,u):z\in B,u<g(z)\}\subset \mathbb{R}^{d}\times \mathbb{R}_{+} $ and
that
\begin{equation*}
\alpha _{g,f}(s)=\int_{\mathbb{R}^{d}}(1-e^{-sf(z)})d\nu _{g}(dz),\quad \beta
_{B,g,f}(s)=\int_{B^{c}}(1-e^{-sf(z)})d\nu _{g}(dz).
\end{equation*}
We have the following result.
\begin{lem}\label{prelemdet}
Let $U_t=t \int_{B^c} f(z) d\nu_{g}(z)$, then $\forall p \geq 1$
\begin{equation}
E(\frac{1}{(N_{t}(1_{B_{g}}f)+U_{t})^{p}})\leq \frac{1}{\Gamma (p)}
\int_{0}^{\infty }s^{p-1}\exp (-t\alpha _{g,f}(s))ds = \frac{1}{\Gamma (p)} I_t^p(f).  \label{predet}
\end{equation}
Suppose moreover that for some $0<\theta \leq \infty$ 
\begin{equation}
\underline{\lim }_{a\rightarrow \infty }\frac{1}{\ln a}\nu _{g}(f\geq 
\frac{1}{a})= \theta ,  \label{A3}
\end{equation}
then for every $t>0$ and $p\geq 1$ such that $ p/t < \theta$
\begin{equation*}
E(\frac{1}{(N_{t}(1_{B_{g}}f)+U_{t})^{p}})<\infty .
\end{equation*}
\end{lem}
Observe that if $\nu(B) < \infty$ then $E \frac{1}{(N_{t}(1_{B_{g}}f)^p}= \infty$
\begin{pf}
By a change of variables we obtain for every $\lambda >0$ 
\begin{equation*}
\lambda ^{-p}\Gamma (p)=\int_{0}^{\infty }s^{p-1}e^{-\lambda s^{}}ds.
\end{equation*}
Taking the expectation in the previous equality with $\lambda
=N_{t}(f1_{B_{g}})+U_{t}$ we obtain 
\begin{equation*}
E(\frac{1}{(N_{t}(f1_{B_{g}})+U_{t})^{p}})=\frac{1}{\Gamma (p)}
\int_{0}^{\infty }s^{p-1}E(\exp (-s(N_{t}(f1_{B_{g}})+U_{t}))ds.
\end{equation*}
Now from Lemma \ref{laplace} we have
\begin{equation*}
E(\exp (-s N_{t}(f1_{B_{g}})))=\exp (-t(\alpha _{g,f}(s)-\beta
_{B,g,f}(s))).
\end{equation*}
Moreover, from the definition of $U_t$ one can easily check that  $\exp (-s U_{t})\leq \exp (-t\beta
_{B,g,f}(s))$ and then
\begin{equation*}
E(\exp (-s (N_{t}(f1_{B_{g}})+U_{t}))\leq \exp (-t\alpha _{g,f}(s))
\end{equation*}
this achieves the proof of (\ref{predet}). The second part of the lemma
follows directily from  lemma \ref{prelem}. 
\end{pf}

\subsubsection{The Malliavin covariance matrix}
In this section, we prove that under some additional assumptions on $p$ and $t$,  $E(\det \sigma(F_M))^{-p} \leq C_p$, for the Malliavin covariance matrix
$\sigma(F_M)$  defined in section $3.4$. 

We first remark that from Hypothesis 3.2 ii) the tangent flow of equation $(\ref{eq4})$ is invertible and
that the moments of all order of this inverse are finite. More precisely
we
define $Y_{t}^{M},t\geq 0$ and $\widehat{Y}_{t}^{M},t\geq 0$ as the
matrix solutions of the equations
\begin{eqnarray}
Y_{t}^{M} &=&I+\sum_{k=1}^{J_{t}^M}\nabla _{x}c_{M}(\overline{Z}_{k},\overline{
X}_{T_{k}-}^{M})Y_{T_{k}-}^{M}+\int_{0}^{t}\nabla _{x}g(\overline{X}
_{s}^{M})Y_{s}^{M}ds, \label{flot}\\
\widehat{Y}_{t}^{M} &=&I-\sum_{k=1}^{J_{t}^M}\nabla _{x}c_{M}(I+\nabla
_{x}c_{M})^{-1}(\overline{Z}_{k},\overline{X}_{T_{k}-}^{M})\widehat{Y}
_{T_{k}-}^{M}-\int_{0}^{t}\nabla _{x}g(\overline{X}_{s}^{M})\widehat{Y}
_{s}^{M}ds. \label{flotinv}
\end{eqnarray}
Then $\widehat{Y}_{t}^{M}\times Y_{t}^{M}=I,\forall t\geq 0.$ Moreover we can prove
under hypotheses $3.0$, $3.1$ and $3.2.$ that  $\forall p \geq 1$ 
\begin{equation}
E(\sup_{s\leq t}(\left\Vert \widehat{Y}_{s}^{M}\right\Vert ^{p}+\left\Vert
Y_{s}^{M}\right\Vert ^{p}))\leq K_{p}<\infty  \label{n14}
\end{equation}
where $K_{p}$ is a constant.
\begin{lem}\label{lemsigmaFM}
Assuming hypothesis $3.0$, $3.1$, $3.2$  we have for  $ p \geq 1$,  $t>0$ such that $ 2dp/ t<\theta$
\begin{equation}
E(\frac{1}{(\det \sigma (F_M))^{p}})\leq C_p,  \label{sigmaFM}
\end{equation}
where the constant $C_p$ does not depend on $M$.
\end{lem}

\begin{pf}
 We first give a lower bound for the lowest eigenvalue of the matrix $\sigma(\overline{X}_t^M)$.
\begin{equation*}
\rho _{t }:=\inf_{\left\vert \xi \right\vert =1}\left\langle
\sigma (\overline{X}_t^M)\xi ,\xi \right\rangle =\inf_{\left\vert \xi
\right\vert =1}\sum_{k=1}^{J_{t}^M} \sum_{r=1}^{d}\left\langle D_{k,r}
\overline{X}_{t}^{M},\xi \right\rangle ^{2}.
\end{equation*}
But from equation $(\ref{eq4})$ we have
$$
D_{k,r} \overline{X}_t^M=\sum_{k'=1}^{J_t^M} \nabla_z c_M(\overline{Z}_{k'}, \overline{X}_{T_{k'}^-}^M) D_{k,r} \overline{Z}_{k'}
+\sum_{k'=1}^{J_t^M} \nabla_x c_M(\overline{Z}_{k'}, \overline{X}_{T_{k'}^-}^M) D_{k,r} \overline{X}_{T_{k'}^-}^M
+ \int_0^t \nabla_x g(\overline{X}_s^M) D_{k,r} \overline{X}_s^M ds
$$
where $\nabla_z c_M=(\partial_{z_r} c_M^{r'})_{r',r}$ and $\nabla_x c_M=(\partial_{x_r} c_M^{r'})_{r',r}$. Since $ D_{k,r} \overline{Z}_{k'}=0$ for $k \neq k'$ we obtain
\begin{equation*}
D_{k,r}\overline{X}_{t}^{M,r'}=(Y_t^{M}\nabla _{z} c_M(\overline{Z}_{k}, \overline{X}_{T_{k}^-}^M) D_{k,r} \overline{Z}_{k})_{r',r}=
\pi_k (Y_t^{M}\nabla _{z} c_M(\overline{Z}_{k}, \overline{X}_{T_{k}^-}^M) )_{r',r} .
\end{equation*}
We deduce that
\begin{equation*}
\sum_{r=1}^{d} \left\langle D_{k,r}\overline{X}_{t}^{M},\xi \right\rangle
^{2} =\sum_{r=1}^{d} \pi_k^2 \left\langle \partial _{z^{r}} c_M(\overline{Z}_{k}, \overline{X}_{T_{k}^-}^M),(Y_{t}^{M})^{\ast }\xi \right\rangle
^{2},
\end{equation*}
but  recalling that $\pi_k \geq 1_{B_{M-1}}(\overline{Z}_k)$ and $c_M=c$ on $B_{M-1}$ we obtain
$$
\sum_{r=1}^{d}  \left\langle D_{k,r}\overline{X}_{t}^{M},\xi \right\rangle
^{2} 
\geq \sum_{r=1}^{d}1_{B_{M-1}}(\overline{Z}_k) \left\langle \partial _{z^{r}} c(\overline{Z}_{k}, \overline{X}_{T_{k}^-}^M),(Y_{t}^{M})^{\ast }\xi \right\rangle
^{2},
$$
and consequently using hypothesis $3.2. iii)$ 
$$
\rho_t \geq \inf_{\left\vert \xi \right\vert =1}
\sum_{k=1}^{J_t^M} 1_{B_{M-1}}(\overline{Z}_k)  \underline{c}^2(\overline{Z}_{k}) \vert (Y_{t}^{M})^{\ast }\xi \vert^2
\geq \left\Vert \widehat{Y}_{t}^{M}\right\Vert ^{-2} \sum_{k=1}^{J_t^M} 1_{B_{M-1}}(\overline{Z}_k)  \underline{c}^2(\overline{Z}_{k}).
$$
Now since $\sigma(F_M)=\sigma(\overline{X}_t^M)+U_M(t)$ we have
$$
E \left\vert \frac{1}{\det \sigma(F_M)} \right\vert^p \leq E \left\vert \frac{1}{\rho_t+U_M(t)} \right\vert^{dp}\leq E\left(\frac{1+\left\Vert \widehat{Y}_{t}^{M}\right\Vert ^{2}}{\sum_{k=1}^{J_t^M} 1_{B_{M-1}}(\overline{Z}_k)  \underline{c}^2(\overline{Z}_{k})+U_M(t)} \right)^{dp}.
$$
Now observe that the denominator of the last fraction is equal in law to
$$
\sum_{k=1}^{J_t^M} 1_{B_{M-1}}(Z_k)  \underline{c}^2(Z_{k})1_{U_k<\gamma(Z_k, X_{T_k-}^M)}+U_M(t) \geq N_t(1_{B^M_{\underline{\gamma}}} \underline{c}^2) +U_M(t),
$$
with $B^M_{\underline{\gamma}}=\{(z,u); z \in B_{M-1}; 0<u<\underline{\gamma}(z)\}$. Assuming hypothesis $3.2.iii$, we can apply 
lemma \ref{prelemdet} with $f=\underline{c}^2$ and $d\nu(z)=\underline{\gamma}(z) d \mu(z)$. This gives for $ p' \geq 1$ such that $p'/t< \theta$
$$
E\left(\frac{1}{ N_t(1_{B^M_{\underline{\gamma}}} \underline{c}^2) +U_M(t)}\right)^{p'} \leq C_{p'}.
$$
Finally since the moments of $\left\Vert \widehat{Y}_{t}^{M}\right\Vert $ are bounded uniformly on $M$ the result of lemma \ref{lemsigmaFM}
follows from  Cauchy-Schwarz inequality.

\end{pf}

\section{References}

[B] V. Bally: An elementary introduction to Malliavin calculus. Preprint No
4718, INRIA February 2003.

[B.B.M] V. Bally, M-P. Bavouzet and M. Messaud: Integration by parts formula
for locally smooth laws and applications to sensitivity computations. Annals
of Applied Probability 2007, Vol. 17, No. 1, 33-66.

[B.F] V.Bally and N.Fournier: regularization properties of the 2D homogeneous
Boltzmann equation without cutoff. Preprint.

[Ba.M] M-P. Bavouzet and M. Messaud: Computation of Greeks using Malliavin
calculus in jump type market models. Electronic Journal of Probability, 11/
276-300, 2006.

[Bi] J.M.\ Bismut: Calcul des variations stochastiques et processus de
sauts. Z.\ Wahrsch. Verw. Gebite, No 2, 147-235, 1983.

[B.G.J] K.\ Bichteler, J.B. Gravereaux and J.\ Jacod: Malliavin calculus for
processes with jumps. Gordon and Breach, 1987.

[Bou] N.\ Bouleau: Error calculus for finance and Physics, the language of
Dirichlet forms. De Gruyer, 2003.

[F.1] N.\ Fournier: Jumping SDE's: Absolute continuity using monotonicity.
SPA, 98 (2), pp 317-330, 2002.

[F.2] N.\ Fournier: Smoothness of the law of some one-dimensional jumping
SDE's with non constant rate of jump. In preparation.

[F.G] N. Fournier and J.S. Giet: On small particles in
coagulation-fragmentation equations. J. Statist. Phys. 111 (5/6) pp
1299-1329, 2003.

[I.W] N.\ Ikeda and S.\ Watanabe: Stochastic Differential Equations and
Diffusion processes. North- Holland, 1989.

[L] R. L\'{e}andre: R\'{e}gularit\'{e} de processus de sauts d\'{e}g\'{e}n\'{e}r\'{e}s. Ann. Inst. H.\ Poincar\'{e}, Proba. Stat., 21, No 2, 125-146,
1985.

[N] D. Nualart: Malliavin calculus and related topics. Springer Verlag, 1995.

\bigskip

\end{document}